\documentclass{amsart}
\usepackage{amssymb, latexsym,epsfig}

\usepackage[dvips]{color}

        \def\NN{\Bbb{N}}

\def\B1{{\rm\kern.32em\vrule    width.12em       height1.4ex
depth-.05ex\kern-.28em 1}}

\newcommand{\R}{{\mathbb R}}
\newcommand{\C}{{\mathbb C}}

\newcommand{\N}{{\mathbb N}}
\newcommand{\CC}{{\overline \C}}

\newcommand{\K}{{\hat{K}}}
\newcommand{\J}{{\mathcal{J}}}
\newcommand{\G}{{\mathcal{G}}}
\newcommand{\Jmin}{{J_{\min }}}
\newcommand{\Jmax}{{J_{\max }}}

\newcommand{\Hmin}{{H_{\min }}}
\newcommand{\stl}{{(\star \lambda)}}
\newcommand{\la}{{\lambda}}
\newcommand{\La}{{\Lambda}}

      \def\NN{\Bbb{N}}
\def\B1{{\rm\kern.32em\vrule    width.12em       height1.4ex
depth-.05ex\kern-.28em 1}}
\newtheorem{thm}{Theorem}[section]

\newtheorem{corollary}[thm]{Corollary}
\newtheorem{lemma}[thm]{Lemma}

\newtheorem{theorem}[thm]{Theorem}
\newtheorem{proposition}[thm]{Proposition}
\newtheorem{claim}[thm]{Claim}

\theoremstyle{definition}
\newtheorem{definition}[thm]{Definition}
\theoremstyle{remark}
\newtheorem{remark}[thm]{Remark}

\newtheorem{example}[thm]{Example}

\numberwithin{equation}{section}

\begin{document}

\title[Dynamics and Structure of Julia sets of polynomial semigroups]
{Dynamical properties and structure of Julia sets of postcritically
bounded polynomial semigroups}

\author{Rich Stankewitz}
\address[Rich Stankewitz]
    {Department of Mathematical Sciences\\
    Ball State University\\
    Muncie, IN 47306\\
    USA}
    \email{rstankewitz@bsu.edu}

\author{Hiroki Sumi}
\address[Hiroki Sumi]
    {Department of Mathematics\\
    Graduate School of Science\\
    Osaka University\\
    1-1, Machikaneyama, Toyonaka\\
    Osaka, 560-0043, Japan}
    \email{sumi@math.sci.osaka-u.ac.jp}

\date{May 9, 2011. To appear in Trans. Amer. Math. Soc.}

\thanks{ 
The first author was partially supported by the BSU Lilly V grant.
He would also like to thank Osaka University for their hospitality
during his stay there while this work was begun.}


\begin{abstract}
We discuss the dynamic and structural properties of polynomial
semigroups, a natural extension of iteration theory to random (walk)
dynamics, where the semigroup $G$ of complex polynomials (under the
operation of composition of functions) is such that there exists a
bounded set in the plane which contains any finite critical value of
any map $g \in G$.  In general, the Julia set of such a semigroup
$G$ may be disconnected, and each Fatou component of such $G$ is
either simply connected or doubly connected (\cite{Su01,Su9}). In this paper, we show
that for any two distinct Fatou components of certain types (e.g.,
two doubly connected components of the Fatou set), the boundaries
are separated by a Cantor set of quasicircles (with uniform
dilatation) inside the Julia set of $G.$  Important in this theory
is the understanding of various situations which can and cannot
occur with respect to how the Julia sets of the maps $g \in G$ are
distributed within the Julia set of the entire semigroup $G$.  We
give several results in this direction and show how such results are
used to generate (semi) hyperbolic semigroups possessing this
postcritically boundedness condition.
\end{abstract}

\maketitle
\section{Introduction}
The dynamics of iteration of a complex analytic map has been studied
quite deeply and in various contexts, e.g., rational, entire, and
meromorphic maps.  It is natural then to consider the generalization
of this theory to the setting where the map may be changed at each
point of the orbit, exactly as in a random walk. Instead of
repeatedly applying the same map over and over again, one may start
with a family of maps $\{ h_{\lambda}: \lambda \in \Lambda\}$, and
consider the dynamics of any iteratively defined \textit{composition
sequence} of maps, that is, any sequence $h_{\lambda_n} \circ \dots
\circ h_{\lambda_1}$ where each $\lambda_k \in \Lambda$. Assigning
probabilities to the choice of map at each stage is the setting for
research of random dynamics (see~\cite{FS, Bruck, Bu1, Bu2, BBR,
Su8, Su02, Su03, Su05, Su9} for previous work related to such dynamics). In this
paper, however, we will concern ourselves with questions of dynamic
stability, not just along such composition sequences one at a time,
but rather we will study when such stability exists no matter which
composition sequence is chosen. Restricting one's attention to the
case where all $h_\lambda$ are rational, one is lead to study the
dynamics of rational semigroups.

A {\bf rational semigroup} is a semigroup generated by non-constant
rational maps on the Riemann sphere $\CC$ with the semigroup
operation being the composition of maps.  We denote by $\langle
h_{\lambda}: \lambda \in \Lambda \rangle $ the rational semigroup
generated by the family of maps $\{ h_{\lambda}: \lambda \in
\Lambda\}.$ A {\bf polynomial semigroup} is a semigroup generated by
non-constant polynomial maps.  Research on the dynamics of rational
semigroups was initiated by A.~Hinkkanen and G.J.~Martin
in~\cite{HM1}, who were interested in the role of the dynamics of
polynomial semigroups while studying various one-complex-dimensional
moduli spaces for discrete groups.  Also, F.~Ren, Z.~Gong, and
W.~Zhou studied such semigroups from the perspective of random
dynamical systems (see~\cite{ZR,GR}). Note that there is a strong
connection between the study of dynamics of rational semigroups and
that of random complex dynamics (see \cite{Su01,Su02,Su03,Su05,Su9}).
For example,
for a random dynamical system generated by a family of polynomials,
we consider the function $T_{\infty }:\overline{\Bbb{C}}\rightarrow [0,1]$
of the probability of tending to $\infty $. It turns out that under certain conditions,
the function $T_{\infty }$ is continuous on $\overline{\Bbb{C}}$ and
varies only on the Julia set $J$ of the associated polynomial semigroup $G$,
and $J$ is a thin fractal set. Moreover, $T_{\infty }$ respects
the surrounding order (see Definition~\ref{def-surr}).
The function $T_{\infty }$ is a complex analogue of
the devil's staircase or Lebesgue's singular functions.
For the detail of these results, see the second author's works \cite{Su05, Su8, Su11}.

As is well known, the iteration of polynomial maps $f_c(z)=z^2+c$
for $c$ in the Mandelbrot set (where the orbit of the sole critical
point $\{f_c^n(0)\}$ is bounded in $\C$), provides a rich class of
maps with many interesting properties.  Many of these dynamic and
structural properties are direct consequences the boundedness of the
critical orbit. It is then natural to look at the more general
situation of polynomial semigroups with bounded postcritical set. We
discuss the dynamics of such polynomial semigroups as well the
structure of their Julia sets. For some properties of polynomial
semigroups with bounded finite postcritical set see also
\cite{Su6,Su5,Su01,Su02,Su03,Su04, Su9}.
This paper is a continuation of the program initiated in~\cite{Su01,Su02,Su03, Su9}.
One may also see the paper~\cite{SS-AMC}, which is based
on conference talks by the authors, where preliminary (weaker) versions
of a portion of the results contained herein were announced and proven.

\begin{definition}
Let $G$ be a rational semigroup. We set
\[ F(G) = \{ z\in \CC \mid G \mbox{ is normal in a neighborhood of  $z$} \} \textrm{ and }
\ J(G)  = \CC \setminus  F(G) .\]We call \(  F(G)\) the {\bf Fatou
set}  of  $G$ and \( J(G)\)  the {\bf Julia set} of $G$. The Fatou
set and Julia set of the semigroup generated by a single map $g$ is
denoted by $F(g)$ and $J(g)$, respectively.
\end{definition}

We quote the following results from~\cite{HM1}. The Fatou set $F(G)$
is \textbf{forward invariant} under each element of $G$, i.e.,
$g(F(G)) \subset F(G)$ for all $g \in G$, and thus $J(G)$ is
\textbf{backward invariant} under each element of $G$, i.e.,
$g^{-1}(J(G)) \subset J(G)$ for all $g \in G$.  Furthermore, when the cardinality $\#J(G)$ is three or more,
$J(G)$ is the smallest closed subset of $\CC$ which contains three or more
points and is backward invariant. Letting the \textbf{backward
orbit} of $z$ be denoted by $G^{-1}(z) =\cup _{g\in G} g^{-1}(\{z\})
$, we have that $J(G)=\overline{G^{-1}(z)}$ for any $z \in J(G)$
whose backward orbit contains three or more points.

We should take a moment to note that the sets $F(G)$ and $J(G)$ are,
however, not necessarily completely invariant under the elements of
$G$. This is in contrast to the case of \textit{iteration} dynamics,
i.e., the dynamics of semigroups generated by a single rational
function. For a treatment of alternatively defined
\textit{completely} invariant Julia sets of rational semigroups the
reader is referred to~\cite{RS,RS1,RS2, SSS}.

 Although the Julia set of a rational semigroup $G$ may not be completely invariant,
 $J(G)$ has an interesting property. Namely, if $G$ is generated by a compact family $\{h_\lambda: \lambda \in \Lambda\}$ of
 rational maps, then $J(G)=\bigcup _{\lambda \in \Lambda }h_\lambda^{-1}(J(G)).$
 This property is called the {\bf backward self-similarity}. In particular,
 if $G=\langle h_{1},\ldots ,h_{m}\rangle $, then
 $J(G)=\bigcup _{j=1}^{m}h_{j}^{-1}(J(G)).$
 From this property, the dynamics of rational semigroups can be regarded as
``backward iterated function systems,'' and in the study of
rational semigroups,  we sometimes borrow and further develop techniques from iterated function systems and fractal geometry.
For these things, see the second author's works \cite{Su1} -- \cite{Su9} and \cite{SU1,SU3}.

Note that $J(G)$ contains the Julia set of each element of $G$.
Moreover, the following critically important result first due to Hinkkanen
and Martin holds (see also~\cite{RS-Bobfest}).

\begin{theorem}[\cite{HM1}, Corollary 3.1] \label{HMJ}
        For rational semigroups $G$ with $\sharp (J(G))\geq 3$,
         we have $$\displaystyle{J(G) = \overline{\bigcup_{f \in G} J(f)}}.$$
\end{theorem}

\begin{remark}
Theorem~\ref{HMJ} can be used to easily show that  $F(\langle
h_\lambda: \lambda \in \Lambda \rangle)$ is \textit{precisely} the
set of $z \in \CC$ which has a neighborhood on which every
composition sequence generated by $\{ h_\lambda: \lambda \in \Lambda
\}$ is normal, i.e., has stable dynamics (see~\cite{ZR}).
\end{remark}

In what follows we employ the following notation.  The {\bf forward
orbit} of $z$ is given by $G(z) =\cup _{g\in G} g(\{z\}) $.  For any
subset $A$ of $\CC ,\ $ we set $G^{-1}(A)=\cup _{g\in G} g^{-1}(A).$
For any polynomial $g$, we denote the \textbf{filled-in Julia set}
of $g$ by $K(g):=\{ z\in \C \mid \cup _{n\in \NN }\, g^{n}(\{z\})
\mbox{ is bounded in }\C \}$.  We note that
for a polynomial $g$ with $\deg (g)\geq 2$,
$J(g)=\partial K(g)$ and
$K(g)$ is the polynomial hull of $J(g)$.  The appropriate
extension (to our situation with polynomial semigroups) of the
concept of the filled-in Julia set is as follows.  (See~\cite{HM1,
Bo-basin} for other kinds of filled-in Julia sets.)

\begin{definition} For a polynomial semigroup $G$,\ we set
$$ \hat{K}(G):=\{ z\in \C
\mid G(z)\mbox{ is bounded in }\C \},$$ and call $\K(G)$ the
\textbf{smallest filled-in Julia set}.
\end{definition}

\begin{remark}
We note that for all $g \in G$, we have $\K(G) \subset K(g)$ and
$g(\K(G)) \subset \K(G)$.
\end{remark}

\begin{definition}
The {\bf postcritical set} of a rational semigroup $G$ is defined by
$$ P(G)=
\overline{\bigcup _{g\in G}\{ \mbox{all critical values of }g:\CC
\rightarrow \CC \} } \ (\subset \CC ).$$
\end{definition}
We say that $G$ is \textbf{hyperbolic} if $P(G) \subset F(G)$ and we
say that $G$ is \textbf{subhyperbolic} if both $\#\{P(G)\cap J(G)\}<
+\infty$ and $P(G) \cap F(G)$ is a compact subset of $F(G).$   For research on
(semi-)hyperbolicity and Hausdorff dimension of Julia sets of
rational semigroups see~\cite{Su1,Su3,Su2,Su7,Su4,Su01,Su02,Su03,SU1, Su9}.

\begin{remark}
It is clear that if rational semigroup $G$ is hyperbolic, then each
$g \in G$ is hyperbolic.  However, the converse is not true.  See
Remark~\ref{ghyp-not-Ghyp}.
\end{remark}

\begin{definition}
The \textbf{planar postcritical set} (or, the finite postcritical
set) of a polynomial semigroup $G$ is defined by
$$P^{\ast }(G)=P(G)\setminus \{ \infty \}.$$
We say that a polynomial semigroup $G$ is \textbf{postcritically
bounded} if $P^*(G)$ is bounded in $\C$.
\end{definition}

\begin{definition}
Let ${\G} $ be the set of all polynomial semigroups $G$ with the
following properties:
\begin{itemize}
\item  each element of $G$ is of degree at least two, and
\item  $P^*(G)$ is bounded in $\C,$ i.e., $G$ is postcritically
bounded.
\end{itemize}
Furthermore, we set ${\G}_{con}= \{ G\in {\G}\mid J(G)\mbox{ is
connected}\} $ and ${\G}_{dis}= \{ G\in {\G}\mid J(G)\mbox{ is
disconnected}\}.$
\end{definition}

\begin{remark}
If $G=\langle h_\lambda: \lambda \in \Lambda \rangle$, then
$P(G)=\overline{\cup_{\lambda \in \Lambda} \cup_{z \in
CV(h_\lambda)} (G(z) \cup \{z\}) }$ where $CV(h)$ denotes the
critical values of $h$.  From this one may, in the finitely
generated case, use a computer to see if $G \in \G$ much in the same
way as one verifies the boundedness of the critical orbit for the
maps $f_c(z)=z^2+c$.  The freely available software~\cite{RS-JuliaSoftware} can be used for this purpose.
\end{remark}

\begin{remark}\label{pinK}
Since $P(G)$ is forward invariant under $G$, we see that $G \in \G$
implies $P^*(G) \subset \K(G)$, and thus $P^*(G) \subset K(g)$ for
all $g \in G$.
\end{remark}

\begin{remark}
For a polynomial $g$ of degree two or more, it is well known that
$\langle g \rangle \in \G$ if and only if $J(g)$ is connected
(see~\cite{Be}, Theorem 9.5.1). Hence, for any $g \in G \in \G$, we
have that $J(g)$ is connected. We note, however, that the analogous
result for polynomial semigroups does not hold as there are many
examples where $G \in \G$, but $J(G)$ is not connected
(see~\cite{SY, Su01,Su02, Su03,Su9}).

 See also \cite{Su04} for an analysis of the
number of connected components of $J(G)$ involving the inverse limit
of the spaces of connected components of the realizations of the nerves of finite
coverings $\mathcal{U}$ of $J(G)$, where $\mathcal{U}$ consists of
backward images of $J(G)$ under finite word maps in $G$.
In fact, the number of connected components of the Julia set of a finitely generated rational semigroup
is deeply related to a new kind of cohomology (so called the ``interaction cohomology''),
which has been introduced by the second author of this paper. Using this cohomology,
one can also investigate the number of connected components of the Fatou set of a finitely generated
rational semigroup.
\end{remark}

The aim of this paper is to investigate what can be said about the
structure of the Julia sets and the dynamics of semigroups $G \in
\G$?  We begin by examining the structure of the Julia set and note
that a natural order (that is respected by the backward action of
the maps in $G$) can be placed on the components of $J(G)$, which
then leads to implications on the connectedness of Fatou
components.\\

\noindent \textbf{Notation:} For a polynomial semigroup $G \in \G$,\
we denote by ${\J}={\J}_G$ the set of all connected components of
$J(G)$ which do not include $\infty$.

\begin{definition}\label{def-surr}
We place a partial order on the space of all non-empty connected
sets in $\C$ as follows.  For connected sets $K_{1}$ and $K_{2}$ in
$\C ,\ $  ``$K_{1}\leq K_{2}$'' indicates that $K_{1}=K_{2}$ or
$K_{1}$ is included in a bounded component of $\C \setminus K_{2}.$
Also, ``$K_{1}<K_{2}$'' indicates $K_{1}\leq K_{2}$ and $K_{1}\neq
K_{2}.$  We call $\leq$ the {\it{surrounding}} order and read $K_1 <
K_2$ as ``$K_1$ is surrounded by $K_2$".
\end{definition}

\noindent \textbf{Convention:} When a set $K_1$ is contained in the
unbounded component of $\CC \setminus K_2$ we say that $K_1$ is
``outside" $K_2$.

\begin{theorem}[\cite{Su01, Su9}]\label{mainth1}
Let $G\in {\G}$ (possibly infinitely generated). Then
\begin{enumerate}
\item \label{mainth1-1}
$({\J},\ \leq )$ is totally ordered.
\item \label{mainth1-2}
Each connected component of $F(G)$ is either simply or doubly
connected.
\item \label{mainth1-3}
For any $g\in G$ and any connected component $J$ of $J(G)$,\ we have
that $g^{-1}(J)$ is connected. Let $g^{\ast }(J)$ be the connected
component of $J(G)$ containing $g^{-1}(J).$ If $J\in {\J}$, then
$g^{\ast }(J)\in {\J}.$ If $J_{1},J_{2}\in {\J} $ and $J_{1}\leq
J_{2},\ $ then both $g^{-1}(J_{1})\leq g^{-1}(J_{2})$ and $g^{\ast
}(J_{1})\leq g^{\ast }(J_{2}).$
\end{enumerate}
\end{theorem}

\begin{figure}[htbp]
\includegraphics[width=2.5cm,width=2.5cm]{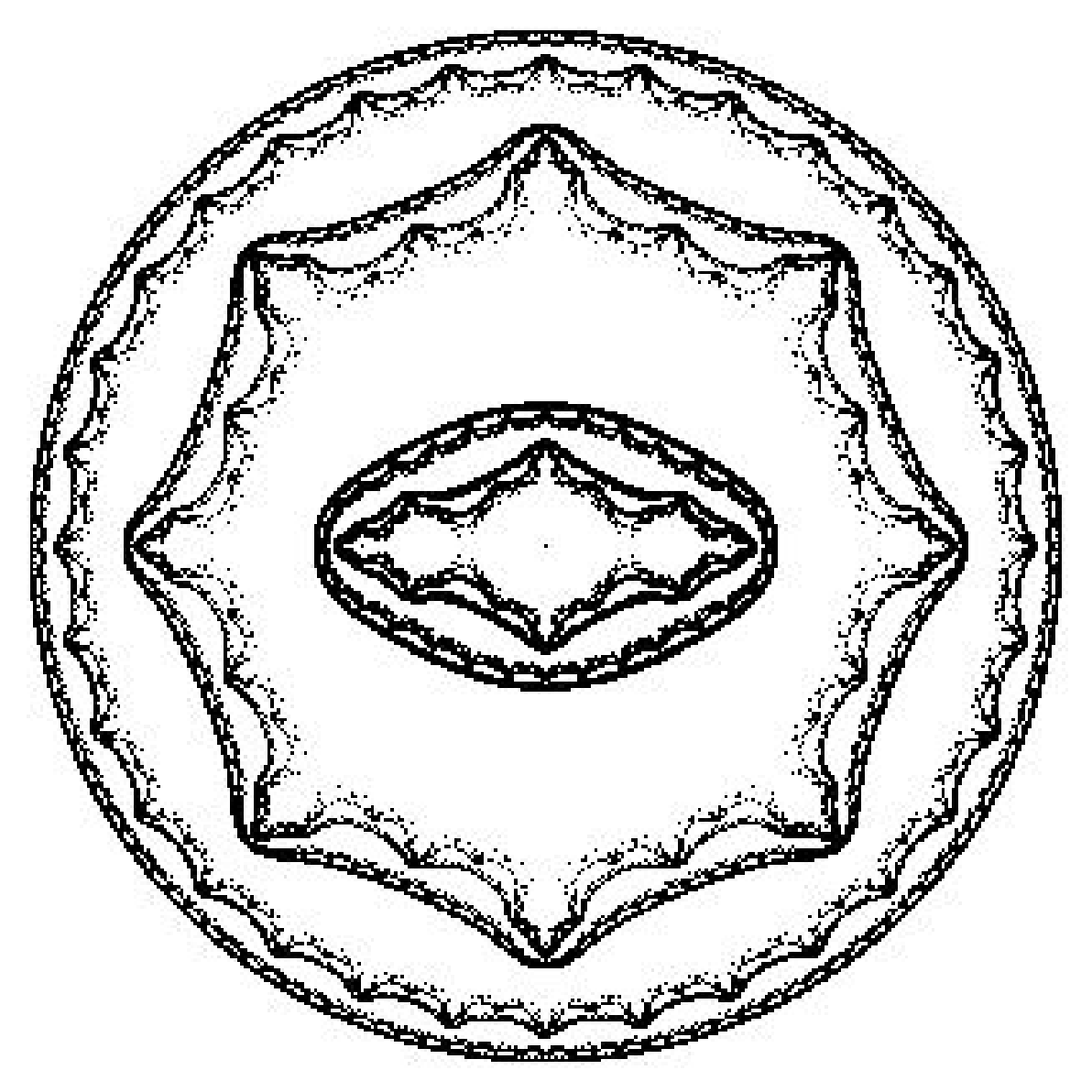}
\caption{The Julia set of $G=\langle h_{1}, h_{2}\rangle $, where we set 
$g_{1}(z):=$ $z^{2}-1, g_{2}(z):=z^{2}/4, h_{1}:=g_{1}^{2}, h_{2}:=
g_{2}^{2}.$
$G\in \mathcal{G}_{dis}$ and $G$ is hyperbolic.\ }
\label{fig:dcjulia}
\end{figure}

\begin{remark}
We note that under the hypothesis of the above theorem $J_1<J_2$ for
$J_1, J_2 \in \J$ does not necessarily imply $g^*(J_1)<g^*(J_2)$, but only that
$g^*(J_1) \leq g^*(J_2)$ as can be seen in
Example~\ref{no-lambda-star}. The Julia set of a $G\in \mathcal{G}_{dis}$ is shown in Figure~\ref{fig:dcjulia}.
\end{remark}

We now present the main results of this paper, first giving
some notation that will be needed to state our
result on the existence of quasicircles in $J(G)$.\\

\noindent \textbf{Notation:}  Given polynomials $\alpha_1$ and
$\alpha_2$, we set $\Sigma_2 = \{x=(\gamma_1, \gamma_2, \dots):
\gamma_k \in \{\alpha_1, \alpha_2\}\}$.  Then, for any $x=(\gamma_1,
\gamma_2, \dots) \in \Sigma_2$, we set $J_x$ equal to the set of
points $z\in \CC $ where the sequence of functions $\{
\gamma_{n}\circ \cdots \circ \gamma_{1}\} _{n\in \NN }$ is not
normal. This is sometimes called the \textit{Julia set along the
trajectory} (sequence) $x \in \Sigma _{2}.$  See~\cite{Su3,Su2,Su4,
Su02,Su03,Su9} for much more on such \textit{fiberwise dynamics}.

\begin{theorem}\label{cantor}
Let $G\in {\mathcal G}_{dis}$ and let $A$ and $B$ be disjoint
subsets of $\CC .$ Suppose that $A$ is a doubly connected component
of $F(G)$ and $B$ satisfies one of the following conditions:
\begin{itemize}
\item $B$ is a doubly connected component of $F(G)$,
\item $B$ is the connected component of $F(G)$ with $\infty \in
B$,
\item $B=\hat{K}(G)$.
\end{itemize}
Then $\partial A\cap \partial B=\emptyset.$  Furthermore,
$\overline{A}$ and $\overline{B}$ are separated by a Cantor family
of quasicircles with uniform dilatation which all lie in $J(G).$
More precisely,\ there exist two elements $\alpha _{1}, \alpha
_{2}\in G$ satisfying all of the following.
\begin{enumerate}
\item There exists a non-empty open set $U$ in $\CC $
with $\alpha _{1}^{-1}(\overline{U})\cap \alpha
_{2}^{-1}(\overline{U})=\emptyset $ and $\alpha
_{1}^{-1}(\overline{U}) \cup \alpha _{2}^{-1}(\overline{U})\subset
U.$
\item $H=\langle \alpha _{1},\alpha _{2}\rangle $
is hyperbolic.
\item
Letting $\Sigma _{2}$ denote the sequence space associated with $\{
\alpha _{1},\alpha _{2}\}$, we have
\begin{enumerate}
     \item $J(H)=\bigcup _{x\in \Sigma _{2}}
     J_{x}$ (disjoint union),\
     \item
     for any component $J$ of $J(H)$ there exists a unique element
     $x \in \Sigma _{2}$ with
     $J=J_{x}$, and
     \item there exists a constant
     $K\geq 1$ such that
     any component $J$ of $J(H)$ is a $K$-quasicircle.

     \end{enumerate}

\item $\{ J_{x}\}
_{x\in \Sigma _{2}}$ is totally ordered with $\leq $,\ consisting of
mutually disjoint subsets of $J(H).$

\item For each $x\in \Sigma _{2}$,\  the set $J_{x}$
separates $\overline{A}$ from $\overline{B}.$
\end{enumerate}
\end{theorem}

\begin{remark}
It should be noted that in the above theorem, the quasicircles
$J_{x}$ are all disjoint components of $J(H)$, but may all lie in
the same component of $J(G)$.  See the proof of
Theorem~\ref{num-J-comp=k}, where a semigroup is constructed such
that there exist only a finite number of components of the Julia
set.
\end{remark}
\begin{remark}
\label{r:nonqc}
There are many hyperbolic polynomial semigroups $G=\langle \alpha _{1},\alpha _{2}\rangle \in \mathcal{G}_{dis}$
such that for a generic $x\in \Sigma _{2}$, the fiberwise Julia set $J_{x}$ is a Jordan curve but not
quasicircle, the unbounded component of $\overline{\Bbb{C}}\setminus J_{x}$ is a John domain, and
the bounded component of $\overline{\Bbb{C}}\setminus J_{x}$ is not a John domain (see \cite{Su03,Su9,Su02}).
See Figure~\ref{fig:dcjulia}.
This phenomenon does not occur in the usual iteration dynamics of a single polynomial.
\end{remark}
\begin{example}\label{cantor-ex}
We give an example of a semigroup $G \in \G$ such that $J(G)$ is a
``Cantor set of round circles".  Let $f_1(z)=az^k$ and $f_2(z)=bz^j$
for some positive integers $k, j \ge 2$. Then, for $|a|^{k-1} \neq
|b|^{j-1}$, the sets $J(f_1)$ and $J(f_2)$ are disjoint circles
centered at the origin. Let $A$ denote the closed annulus between
$J(f_1)$ and $J(f_2)$. For positive integers $m_1$ and $m_2$ each
greater than or equal to 2 (if $k$ and $j$ are not both equal to 2
then $m_1=m_2=1$ will also suffice), we see that the iterates
$g_1=f_1^{m_1}$ and $g_2=f_2^{m_2}$ will yield $A_1=g_1^{-1}(A) \cup
g_2^{-1}(A) \subset A$ where $g_1^{-1}(A) \cap g_2^{-1}(A)
=\emptyset$.  Now iteratively define $A_{n+1}=g_1^{-1}(A_n) \cup
g_2^{-1}(A_n)$ and note that for $G=\langle g_1, g_2 \rangle$ we
have $J(G)=\cap_{n=1}^\infty A_n$, since $J(G)$ is the smallest
closed backward invariant (under each element of $G$) set which
contains three or more points.
\end{example}

For our remaining results we need to note the existence of both a
minimal element and a maximal element in $\J$ and state a few of
their properties.

\begin{theorem}[\cite{Su01, Su9}]\label{Jmin}
Let $G$ be a polynomial semigroup in $\G _{dis}.$ Then there is a
unique element $J_{\min }(G)$ (abbreviated by $J_{\min }$ ) $\in \J$
such that $\Jmin$ meets (and therefore contains) 
$\partial \K (G)$. Also, $\infty \in F(G)$ and there exists a unique element
$J_{\max }(G)$ (abbreviated by $J_{\max}$) $\in \J$ such
that $\Jmax$ meets (and therefore contains) $\partial U_\infty$,
where $U_\infty$ is the simply connected component of $F(G)$ which
contains $\infty$.  Moreover, int$\hat{K}(G)\neq \emptyset .$
Furthermore, we have the following

$\bullet$ $\Jmin \leq J$ for all $J \in \J$,

$\bullet$ $\Jmax \geq J$ for all $J \in \J$,

$\bullet$ $\K(G)$, and therefore $P^*(G)$, is contained in the
polynomial hull of each $J \in \J$.
\end{theorem}

\begin{remark} We see that $\partial \K(G) \subset J(G)$ when $G \in \G$, but, in
general, we do not have $\partial \K(G)=J(G)$, unlike in iteration
theory where $\partial K(g) = J(g)$ for polynomials $g$ of degree
two or more.  In fact, $\partial \K(G)$ might not even equal
$\Jmin(G)$ either (see Example~\ref{Hmin-not}).
\end{remark}

\begin{remark}
When $G \in \G_{con}$ we will use the convention that
$\Jmin=\Jmax=J(G)$ and note that it is still the case that  $\Jmin$
meets $\partial \K$ and $P^*(G)$ is contained in the polynomial hull
of each $J \in \J$. However, it is not necessarily the case that
$\infty \in F(G)$, as exhibited by the example $\langle z^2/n:n \in
\N \rangle$.
\end{remark}

In the proofs of many results concerning postcritically bounded
polynomial semigroups, it is critical to understand the
distribution of the sets $J(g)$ where $g \in G$, especially when $g$
is a generator of $G$.  In particular, it is important to understand
the relationship between such $J(g)$ and the special components
$\Jmin$ and $\Jmax$ of $J(G)$.  In Section~\ref{struct-prop} we
investigate such matters carefully providing several results
including Theorem~\ref{num-J-comp=k} below.

In~\cite{Su01}, it was shown that, for each positive integer $k$,
there exists a semigroup $G\in {\G}_{dis}$ with $2k$ generators such
that $J(G)$ has exactly $k$ components. Furthermore, in~\cite{Su04}
it was shown that any semigroup in ${\G}$ generated by exactly three
elements will have a Julia set with either one component or infinitely many
components (examples where the number of components is one, $\aleph_0$ or uncountable were given). Hence we have the following question: For fixed integer
$k\geq 3$, what is the fewest number of generators that can produce
a semigroup $G\in {\G}_{dis}$ with $\sharp {\J}=k$? The answer to
this question is four as stated in Theorem 1.23 below.

\begin{theorem}\label{num-J-comp=k}
For any $k \in \N$, there exists a 4-generator polynomial semigroup
$H \in \G $ such that $\#\J_H=k$.  Furthermore, $H$ can be
chosen so that no $J \in \J_H\setminus \{\Jmin, \Jmax\}$ meets the
Julia set of any generator of $H$.
\end{theorem}

The next two results, whose proofs depend on understanding the
distribution of the $J(g)$ within $J(G)$, concern the
(semi-)hyperbolicity of polynomial semigroups in $\G$. In
particular, they show how one can build larger (semi-)hyperbolic
polynomial semigroups in $\G$ from smaller ones by including maps
with certain properties.  We first state two definitions.

\begin{definition}
We define Poly $=\{ h:\CC \rightarrow \CC \mid h \mbox{ is a
non-constant polynomial}\} $, endowed with the topology of uniform convergence on
$\CC$ with respect to the spherical metric.
\end{definition}

\begin{remark} \label{Poly}
For use later we note that given integer $d \geq 1$, a sequence
$p_n$ of polynomials of degree $d$ converges to a polynomial $p$ in
Poly if and only if the coefficients converge appropriately and $p$
is of degree $d$.
\end{remark}

\begin{definition} A
rational semigroup $H$ is \textbf{semi-hyperbolic} if for each $z
\in J(H)$ there exists a neighborhood $U$ of $z$ and a number $N \in
\N$ such that for each $g \in H$ we have $\deg(g:V \to U) \leq N$
for each connected component $V$ of $g^{-1}(U)$.
\end{definition}

\begin{theorem}\label{semi-hyp}
Let $H \in \G$, $\Gamma$ be a compact family in {\em Poly}, and let
$G=\langle H, \Gamma \rangle$ be the polynomial semigroup generated
by $H$ and $\Gamma.$
Suppose \\
\indent (1)\, $G\in \G _{dis}$, \\
\indent (2)\, $J(\gamma)\cap J_{\min }(G)=\emptyset$ for each
$\gamma \in \Gamma$,
and \\
\indent (3)\, $H$ is semi-hyperbolic. \\
Then, $G$ is semi-hyperbolic.
\end{theorem}

\begin{remark}
Theorem~\ref{semi-hyp} would not hold if we were to replace both
instances of the word \textit{semi-hyperbolic} with the word
\textit{hyperbolic} (see Example~\ref{nothyp}).  However, with an
additional hypothesis we do get the following result.
\end{remark}

\begin{theorem}\label{hyp}
Let $H \in \G$, $\Gamma$ be a compact family in Poly, and let
$G=\langle H, \Gamma \rangle$ be the polynomial semigroup generated
by $H$ and $\Gamma.$
Suppose \\
\indent (1)\, $G\in {\mathcal G}_{dis}$, \\
\indent (2)\, $J(\gamma)\cap J_{\min }(G)=\emptyset$ for each
$\gamma \in \Gamma$,\\
\indent (3)\, $H$ is hyperbolic, and \\
\indent (4)\, for each $\gamma \in \Gamma$, the critical values of $\gamma$ do not meet $\Jmin(G)$.\\
Then, $G$ is hyperbolic.
\end{theorem}

\begin{remark}
We note that hypothesis (3) can be replaced by the slightly weaker
hypothesis that $P^*(H) \cap J(H) = \emptyset$ since if $\infty \in
J(H) \subset J(G)$, then $J(G)$ is connected by Theorem~\ref{Jmin}
and so hypothesis (1) fails to hold. A similar remark can be made
about hypothesis (3) in Theorem~\ref{semi-hyp}.
\end{remark}

\begin{remark}
Theorems~\ref{semi-hyp} and~\ref{hyp} do not require that $H$ or $G$
be generated by a finite, or even compact, subset of Poly.
\end{remark}

The rest of this paper is organized as follows.  In
Section~\ref{BandT} we give the necessary background and tools
required.  In Section~\ref{pf-q-circles} we give the proof of
Theorem~\ref{cantor}.  In Section~\ref{struct-prop} we provide a
more detailed look at the distribution of $J(g)$ within $J(G)$, in
particular, proving Theorem~\ref{num-J-comp=k}. In
Section~\ref{pf-hyp} we give the proofs of Theorems~\ref{semi-hyp}
and~\ref{hyp} along with Example~\ref{nothyp}.

\section{Background and Tools}\label{BandT}

We first state some notation to be used later.\\

\noindent \textbf{Notation:}  Given any set $A \subset \CC$ we
denote by $\overline{A}$ the closure of $A$ in $\CC$. For $z_0 \in
\C$ and $r, R >0$ we set $B(z_0,r)=\{z \in \C:|z-z_0|<r\}, C(z_0,
r)=\{z \in \C:|z-z_0|=r\}$, and $Ann(z_0;r,R)=\{z \in
\C:r<|z-z_0|<R\}$.  Furthermore, given any set $C \subset \C$ we
denote the $\epsilon-$neighborhood of $C$ by $B(C,
\epsilon)=\cup_{z \in C} B(z,\epsilon)$.  \\

Most often our understanding of the surrounding order $\leq$ given
in Definition~\ref{def-surr} will be applied to \textit{compact}
connected sets in $\C$ and so in this section we state many results
which we will need later. Although not all connected compact sets in
$\C$ are comparable in the surrounding order, we do have the
following two lemmas whose proofs we leave to the reader.

\begin{lemma}\label{surr}
Given two connected compact sets $A$ and $B$ in $\C$ we must have
exactly one of the following:
\begin{enumerate}
\item $A<B$
\item $B<A$
\item $A \cap B \neq \emptyset$
\item $A$ and $B$ are \textbf{outside of each other}, i.e.,
$A$ is a subset of the unbounded component of $\C \setminus B$ and
$B$ is a subset of the unbounded component of $\C \setminus A$.
\end{enumerate}
\end{lemma}

\begin{definition}
For a compact set $A \subset \C$ we define the \textbf{polynomial
hull} $PH(A)$ of $A$ to be the union of $A$ and all bounded
components of $\C \setminus A$.
\end{definition}

\begin{lemma}\label{surr-order}
Let $A$ and $B$ be compact connected subsets of $\C$ such that
$PH(A) \cap PH(B) \neq \emptyset$  Then exactly one of the following
holds:
\begin{enumerate}
\item $A<B$
\item $B<A$
\item $A \cap B \neq \emptyset$.
\end{enumerate}
\end{lemma}

\begin{remark}
We note that for compact connected sets $A$ and $B$ in $\C$, it
follows that $A < B$ if and only if $PH(A) < B$ since the set
$PH(A)$ is also compact and connected.
\end{remark}

\begin{lemma}\label{surr-phull}
Let $g$ be a polynomial of degree at least one and suppose $B
\subset PH(A)$ where $g(B) \subset B$ and $A \subset \C$ is compact.
Then $B \subset PH(g^{-1}(A))$.  In particular, if $g \in G \in \G$
and $P^*(G) \subset PH(A)$ where $A \subset \C$ is compact, then
$P^*(G) \subset PH(g^{-1}(A))$.
\end{lemma}

\begin{proof} Suppose $B \nsubseteq PH(g^{-1}(A))$.  Thus there exists $z_0 \in B$
in the unbounded
component $U$ of $\CC \setminus g^{-1}(A)$.  Let $\gamma$ be a curve
in $U$ connecting $z_0$ to $\infty$.  Then $\Gamma = g \circ \gamma$
is a curve in $\CC \setminus A$ which connects $g(z_0)$ to $\infty$
which shows that $g(z_0) \notin PH(A)$.  Since $B$ is forward
invariant we have that $g(z_0) \in B \setminus PH(A)$ which
contradicts our hypothesis. \end{proof}

\begin{corollary}\label{pc-surr}
Let $f, g \in G \in \G$.  If $A$ is of the form $J \in \J, J(f),
g^{-1}(J)$, or $g^{-1}(J(f))$, then $P^*(G) \subset PH(A)$.
\end{corollary}

\begin{proof} By Theorem~\ref{Jmin} we have $P^*(G) \subset
PH(J)$ for all $J \in \J$. By Remark~\ref{pinK} $P^*(G) \subset
K(f)=PH(J(f))$.  The other cases then follow from
Lemma~\ref{surr-phull}. \end{proof}

\begin{lemma}[\cite{N}]
\label{nadlem} Let $X$ be a compact metric space and let
$f:X\rightarrow X$ be a continuous open map. Let $K$ be a compact
connected subset of $X.$ Then for each connected component $B$ of
$f^{-1}(K)$, we have $f(B)=K.$
\end{lemma}

\begin{lemma}\label{backimage}
Let $g$ be a polynomial with $d=\deg(g) \geq 1$ and let $K \subset
\C$ be a connected compact set such that the unbounded component $U$
of $\CC \setminus K$ contains no critical values of $g$ other than
$\infty$, i.e., the finite critical values of $g$ lie in $PH(K)$.
Then $g^{-1}(K)$ is connected. Further, if $K_1$ is a connected
compact set such that $K<K_1$, then $g^{-1}(K)<g^{-1}(K_1)$.
\end{lemma}

\begin{proof} Set $V=g^{-1}(U)$ and note that $V$ contains no
finite critical points of $g$.  Thus by the Riemann-Hurwitz relation
we have $\chi(V) + \delta_g(V)=d \chi(U)$, where $\chi(\cdot)$
denotes the Euler characteristic and $\delta_g(\cdot)$ is the
deficiency.  Since the hypotheses on $U$ imply $\delta_g(V)=d-1$ and
$\chi(U)=1$, we see that $\chi(V)=1$.  Hence the open and connected
set $V$ is simply connected.

Suppose that $g^{-1}(K)$ is not connected.  Then there exists a
bounded component $V_0$ of $\C \setminus g^{-1}(K)$ which is not
simply connected (see~\cite{Be}, Proposition 5.1.5).  Thus there
exists a Jordan curve $\gamma \subset V_0$ such that the bounded
component $B$ of $\C \setminus \gamma$ contains some component $E$
of $g^{-1}(K)$. Hence $V_0 \cup B$ is open and does not meet $V$. By
Lemma~\ref{nadlem} we have $g(E)=K$. Hence $g(V_0 \cup B) \supset K
\supset \partial U$, which, by the Open Mapping Theorem, implies
$g(V_0 \cup B)$ meets $U$, and thus $V_0 \cup B$ meets $V$.  This
contradiction implies that $V_0$ is simply connected and hence
$g^{-1}(K)$ is connected.

Now suppose $K<K_1$.  Let $U_0$ be the bounded component of $\C
\setminus K_1$ such that $U_0 \supset K$ and let $U_1$ be the
unbounded component of $\C \setminus K_1$.  Hence $g^{-1}(U_0)$ does
not meet $V_1=g^{-1}(U_1)$.  Hence $PH(g^{-1}(K_1)) = \C \setminus
V_1 \supset g^{-1}(U_0) \supset g^{-1}(K)$.  Since $g^{-1}(K_1) \cap
g^{-1}(K) = \emptyset$ (which follows from the fact that $K_1 \cap K
= \emptyset$), we conclude that $g^{-1}(K_1)
> g^{-1}(K)$.
\end{proof}

\begin{corollary}\label{backcor}
Let $g, h \in G \in \G$ and $J \in \J$.  Then $g^{-1}(J)$ and
$g^{-1}(J(h))$ are connected.  Furthermore, $J_1<J_2$ for $J_1, J_2
\in \J$ implies $g^{-1}(J_1)<g^{-1}(J_2)$, and $J(h_1)<J(h_2)$ for
$h_1, h_2 \in G$ implies $g^{-1}(J(h_1))<g^{-1}(J(h_2)).$
\end{corollary}

\begin{proof} Any finite critical value of $g$ must lie in
$P^*(G) \subset PH(J) \cap PH(J(h))$ by Corollary~\ref{pc-surr}. The
result then follows immediately from Lemma~\ref{backimage}.
\end{proof}

\begin{corollary}\label{comp-surr}
Let $f, g \in G \in \G$.  For any two sets $A$ and $B$ of the form
$J \in \J, J(f), g^{-1}(J)$, or $g^{-1}(J(f))$, exactly one of the
following must hold:
\begin{enumerate}
\item $A<B$
\item $B<A$
\item $A \cap B \neq \emptyset$.
\end{enumerate}
\end{corollary}

\begin{proof} This is immediate from
Corollary~\ref{backcor}, Corollary~\ref{pc-surr} and
Lemma~\ref{surr-order}. \end{proof}

The following lemma will allow one to understand the surrounding
order through an imbedding, of sorts, into the real numbers.

\begin{lemma}\label{dist-surr}
Suppose $z_0 \in PH(C_1) \cap PH(C_2)$ where $C_1$ and $C_2$ are
disjoint compact connected sets in $\C$.  Then $dist(z_0, C_1)<
dist(z_0, C_2)$ if and only if $C_1 < C_2$ in the surrounding order.
\end{lemma}

\begin{proof}
First suppose that $d_2=dist(z_0, C_2)>dist(z_0, C_1)$.  Then we
have $B(z_0, d_2) \subset PH(C_2)$.  Since $B(z_0, d_2) \cap C_1
\neq 0$ and $C_1 \cap C_2 = \emptyset$, we must have that the
bounded component of $\CC \setminus C_2$ which contains the
connected set $B(z_0, d_2)$ also contains the connected set $C_1$.
Thus $C_2>C_1$.

Suppose $C_1<C_2.$  Letting $d_1=dist(z_0, C_1)$ we see that
$\overline{B(z_0, d_1)} \subset PH(C_1)<C_2$ implies
$\overline{B(z_0, d_1)}$ must not meet $C_2$, i.e.,
$d_1<dist(z_0,C_2)$.
\end{proof}

\begin{lemma}\label{inf-sup-exist}
Let $\{C_\alpha\}_{\alpha \in \mathcal{A}}$ be a collection of
non-empty compact connected sets in $\C$ that are linearly ordered
by the surrounding order $\leq$.  Suppose $\{C_\beta\}_{\beta \in
\mathcal{B}}$ is a sub-collection of $\{C_\alpha\}_{\alpha \in
\mathcal{A}}$ such that $\overline{\cup_{\beta \in \mathcal{B}}
C_{\beta}} \subset \cup_{\alpha \in \mathcal{A}} C_\alpha$.  Then
both $\inf_{\beta \in \mathcal{B}} C_\beta$ and $\sup_{\beta \in
\mathcal{B}} C_\beta$ exist and are in $\{C_\alpha\}_{\alpha \in
\mathcal{A}}$.
\end{lemma}

\begin{proof}
By compactness and the linear ordering on $\{C_\alpha\}_{\alpha \in
\mathcal{A}}$, one can quickly show that the collection
$\{PH(C_\alpha)\}_{\alpha \in \mathcal{A}}$ satisfies the finite
intersection property. Thus there exists some $z_0 \in \cap_{\alpha
\in \mathcal{A}} PH(C_\alpha)$. For each $\beta \in \mathcal{B}$,
let $r_\beta=dist(z_0,C_\beta)$ and consider $r_0=\inf r_\beta$.

We only need to consider the case where $r_0 < r_\beta$ for all
$\beta \in \mathcal{B}$, since if $r_{\beta_0}=r_0$, then clearly,
by Lemma~\ref{dist-surr}, $C_{\beta_0}=\inf_{\beta \in \mathcal{B}}
C_\beta$. Select a strictly decreasing sequence $r_{\beta_n} \to
r_0$. By Lemma~\ref{dist-surr}, we have that $C_{\beta_1} >
C_{\beta_2} > \dots$.  Let $z_{\beta_n} \in C_{\beta_n}$ be
arbitrary. Without loss of generality we may assume that
$z_{\beta_n} \to a_0 \in \C$. By hypothesis there exists
$C_{\alpha_0}$ which contains $a_0$.  We will now show that
$C_{\alpha_0}=\inf_{\beta \in \mathcal{B}} C_\beta$.

Fixing $\beta \in \mathcal{B}$ and applying Lemma~\ref{dist-surr},
we see that $\{a_0\} < C_\beta$, since the sequence
$(z_{\beta_k})_{k \geq n}$ must lie in $PH(C_{\beta_n})<C_\beta$ for
large $n$ (whenever $r_{\beta_n} < r_\beta$). Thus we must have
$C_{\alpha_0} < C_\beta$ for all $\beta \in \mathcal{B}$. Hence
$C_{\alpha_0}$ is a lower bound for $\{C_\beta\}_{\beta \in
\mathcal{B}}$.  Suppose that $C_{\alpha_1} > C_{\alpha_0}$.  It must
then be the case that $\{a_0\} < C_{\alpha_1}$ and so it follows
that $\{z_{\beta_n}\} < C_{\alpha_1}$ for large $n$.  Thus
$C_{\beta_n} < C_{\alpha_1}$ for large $n$, implying that
$C_{\alpha_1}$ is not a lower bound for $\{C_\beta\}_{\beta \in
\mathcal{B}}$.  We conclude that $C_{\alpha_0}=\inf_{\beta \in
\mathcal{B}} C_\beta$.

The proof that $\sup_{\beta \in \mathcal{B}} C_\beta$ exist in
$\{C_\alpha\}_{\alpha \in \mathcal{A}}$ follows a similar argument
using $\sup r_\beta$ and Lemma~\ref{dist-surr}.  We omit the
details.
\end{proof}

By the proof of the above lemma we see that if
$\overline{\cup_{\beta \in \mathcal{B}} C_{\beta}} =\cup_{\beta \in
\mathcal{B}} C_{\beta}$, then both $\inf_{\beta \in \mathcal{B}}
C_\beta$ and $\sup_{\beta \in \mathcal{B}} C_\beta$ are in
$\{C_\beta\}_{\beta \in \mathcal{B}}$.  Thus we have the following.

\begin{lemma}\label{min-max-exist}
Let $\{C_\beta\}_{\beta \in \mathcal{B}}$ be a collection of compact
connected sets in $\C$ that are linearly ordered by the surrounding
order $\leq$.  If $\overline{\cup_{\beta \in \mathcal{B}} C_{\beta}}
=\cup_{\beta \in \mathcal{B}} C_{\beta}$, then we can conclude that
both $\min_{\beta \in \mathcal{B}} C_\beta$ and $\max_{\beta \in
\mathcal{B}} C_\beta$ exist.
\end{lemma}

\begin{lemma}\label{pullback-order}
Let $f \in G \in \G$.  Let $K$ be a connected compact set in $\C$
such that $PH(K) \supset P^*(f)$.

\noindent (a) Let $J(f) > K$.  Then $J(f) > f^{-1}(K)$.  Also,
$f^{-1}(K)> K$ or $f^{-1}(K) \cap K \neq \emptyset.$

\noindent (b) Let $J(f) < K$.  Then $J(f)<f^{-1}(K)$.  Also,
$f^{-1}(K)<K$ or $f^{-1}(K) \cap K \neq \emptyset.$
\end{lemma}

\begin{proof}
We now prove (a).  We first note that $J(f)=f^{-1}(J(f))
> f^{-1}(K)$ follows immediately from Lemma~\ref{backimage}. Since
$P^*(f) \subset PH(K)<J(f)$ we see that $f$ cannot have a Siegel
disk or parabolic fixed point. Hence, $f$ must have a finite
attracting fixed point $z_0$. Furthermore, since $PH(K)$ is
connected and $P^*(f) \subset PH(K)<J(f)$, it is clear that there
can be only one attracting fixed point for $f$ and $K$ must lie in
the immediate attracting basin $A_f(z_0)$. Since $PH(f^{-1}(K))
\supset P^*(f)$ by Lemma~\ref{surr-phull}, we see that $f^{-1}(K)$
also lies in $A_f(z_0)$. Hence $A_f(z_0)$ must be completely
invariant under $f$. This implies $F(f)$ has only two components
$A_f(\infty)$ and $A_f(z_0)$, each which are simply connected
(see~\cite{Be}, Theorem 5.6.1).

Letting $\varphi:A_f(z_0) \to B(0,1)$ be the Riemann map such that
$z_0 \mapsto 0$, then one may apply Schwarz's Lemma to the degree
greater than or equal to two (finite Blaschke product) map
$B=\varphi \circ f \circ \varphi^{-1}$ to show that any point mapped
to a point of maximum modulus of $\varphi(K)$ must lie outside of
$\varphi(K)$. Thus either $f^{-1}(K)> K$ or $f^{-1}(K) \cap K \neq
\emptyset.$

Part (b) is proved more easily than (a) since it is already known
that $A_f(\infty)$ is simply connected.  Then one can similarly
examine the Riemann map from $A_f(\infty)$ to $B(0,1)$ such that
$\infty \mapsto 0$.
\end{proof}

We note that Theorem~\ref{Jmin} along with the proof of part (a)
above, with $K=\Jmin$, proves the following (which has been already shown
in~\cite{Su01, Su9}).

\begin{lemma}\label{attr-basin}
Let $f \in G \in \G$ be such that $J(f) > \Jmin$, i.e., $J(f) \cap
\Jmin = \emptyset$.  Then $f$ has an attracting fixed point $z_0 \in
\C$ and $F(f)$ consists of just two simply connected immediate
attracting basins $A_f(\infty)$ and $A_f(z_0)$.
\end{lemma}

We note that the maps $f=h_1^2$ and $g=h_2^2$ where $h_1(z)=z^2-1$ and $h_2(z)=z^2/4$
generate $G=\langle f, g \rangle \in \G_{dis}$ where $f$ has two finite attracting fixed points
(see \cite{Su03,Su9}).
Thus we see that the condition $J(f) \cap
\Jmin = \emptyset$ in the above lemma is indeed necessary.
We also note that for this $G$, the phenomena in Remark~\ref{r:nonqc} holds.
See Figure~\ref{fig:dcjulia}.
\begin{lemma}\label{connect-to-Jf}
Let $f \in G \in \G$.  Let $K$ be a connected set in $J(G)$
containing three or more points such that $f^{-n}(K)$ is also
connected for each $n \in \N$. If $f^{-1}(K) \cap K \neq \emptyset$,
then $J(f)$ and $K$ are contained in the same component $J \in \J$.
\end{lemma}

\begin{proof}
The lemma follows from the fact that the connected set
$\overline{\cup_{n =0}^\infty f^{-n}(K)}$ in $J(G)$ must meet
$J(f)$.
\end{proof}

We now present a general topological lemma that will be used to
justify a corollary which will be needed later.
\begin{lemma}\label{generic-isol}
Let $\{C_\alpha\}_{\alpha \in \mathcal{A}}$ be a collection of
compact connected sets in $\CC$. Let $\epsilon>0$ and let $C$ be any
connected component of $\overline{\cup_{\alpha \in \mathcal{A}}
C_\alpha}$. Then there exists $\alpha \in \mathcal{A}$ such that
$C_\alpha \subset B(C, \epsilon)$.
\end{lemma}

\begin{proof}
Choose any $z \in C$ and let $\alpha_n \in \mathcal{A}$ be such that
$dist(z, C_{\alpha_n}) \to 0$.  By compactness in the topology
generated by the Hausdorff metric on the space of non-empty compact
subsets of $\CC$, we then may conclude (by passing to subsequence if
necessary) that $C_{\alpha_n} \to K$ for some non-empty connected
compact set $K$, which therefore must contain $z$ and hence be
contained in $C$. Thus for large $n$ we have $B(C,\epsilon) \supset
B(K,\epsilon) \supset C_{\alpha_n}$.
\end{proof}

Using the fact that $J(g)$ is connected whenever $g \in G \in \G$ we
clearly obtain the following slight generalization of Lemma 4.2
in~\cite{Su01}.
\begin{corollary}\label{isol}
Let $\{g_\la\}_{\la \in \La} \subset G \in \G$.  Let $\epsilon>0$
and let $C$ be any connected component of $\overline{\cup_{\la \in
\La} J(g_\la)}$.  Then there exists $\la \in \La$ such that
$J(g_\la) \subset B(C, \epsilon)$.

In particular, we apply Theorem~\ref{HMJ} to obtain that if
$\{g_\la\}_{\la \in \La} = G \in \G$ and $J \in \J$, then for every
$\epsilon>0$ there exists $g \in G$ such that $J(g) \subset
B(J,\epsilon)$.
\end{corollary}

\section{Proof of Theorem~\ref{cantor}}\label{pf-q-circles}
We first present a definition and a lemma that will assist in the
proof of Theorem~\ref{cantor}.

\begin{definition}
For compact connected sets $K_1$ and $K_2$ in $\C$ such that
$K_1<K_2$ we define $Ann(K_1,K_2)=U \setminus PH(K_1)$ where $U$ is
the bounded component of $\C \setminus K_2$ which contains $K_1$.
Thus $Ann(K_1,K_2)$ is the open doubly connected region ``between"
$K_1$ and $K_2$.
\end{definition}

\begin{remark}
For any compact connected set $A \subset Ann(K_1,K_2)$ we
immediately see that $A<K_2$ and, by Lemma~\ref{surr}, either $K_1$
and $A$ are outside of each other or $K_1<A$.
\end{remark}

\begin{lemma}\label{J(k)}
Let $f, g \in G \in \G$ be such that $J(f)$ and $J(g)$ lie in
different components of $J(G)$ with $J(f)<J(g)$.  Then for any fixed
$n, m \in \N$ there exists $h, k \in G$ such that $f^{-(n+1)}(J(g))
< J(h)< f^{-n}(J(g))$ and $g^{-m}(J(f)) < J(k)< g^{-(m+1)}(J(f))$.
\end{lemma}

\begin{proof} Corollary~\ref{pc-surr}, Lemma~\ref{pullback-order}(a),
and Lemma~\ref{connect-to-Jf} show that $g^{-1}(J(f))>J(f)$. Set
$X=g^{-1}(J(f)), A=g^{-m}(J(f))$ and $B=g^{-(m+1)}(J(f))$ and note
that $J(f)<A<B$ by Lemma~\ref{backimage}. Keeping
Lemma~\ref{attr-basin} in mind, we may choose $\ell \in \N$ large
enough so that $f^{-\ell}(B) \subset Ann(J(f),X)$. Then
$g^{-m}(f^{-\ell}(\overline{Ann(A,B)})) \subset
g^{-m}(Ann(J(f),X))\subset Ann(A,B)\subset \overline{Ann(A,B)}$
which implies that $k=f^{\ell}\circ g^m \in G$ is such that $J(k)
\subset Ann(A,B).$  Since by Corollary~\ref{comp-surr} we must have
either $J(k)<A$ or $A<J(k)$, we see by construction that $A<J(k)$
must hold.

The other result is proved similarly. \end{proof}

We will require the following result which was proved via
\textit{fiberwise} quasiconformal surgery by the second author.

\begin{proposition}[\cite{Su03}, Proposition 2.25] \label{qc-surg}
Let $G=\langle \alpha_1, \alpha_2 \rangle \in \G$ be hyperbolic such
that $P^*(G)$ is contained in a single component of $int (\K(G))$.
Then there exists $K \geq 1$ such that for all sequences $x \in
\Sigma_2$, the set $J_x$ is a $K$-quasicircle.
\end{proposition}

\begin{remark}
Under the hypotheses above we know that for each $g \in G$, the set
$J(g)$ is a quasicircle (see~\cite{CG}, p.~102).  But the above
result shows much more as it shows that the Julia sets along
\textit{sequences} are also all quasicircles, and that all such
quasicircles have uniform dilations.
\end{remark}

We now can present the proof of Theorem~\ref{cantor}.

\begin{proof}[Proof of Theorem~\ref{cantor}] We first give a proof
in the case that $A$ and $B$ are doubly connected components of
$F(G)$.  Since the doubly connected components of $F(G)$ are
linearly ordered by $\leq$, we may assume without loss of generality
that $A < B$.

Let $\gamma_A$ be a non-trivial curve in $A$ (i.e., $\gamma_A$
separates the components of $\C \setminus A$) and let $\gamma_B$ be
a non-trivial curve in $B$.  Since $J(G)=\overline{\cup_{g \in G}
J(g)}$, the bounded component of $\C \setminus \gamma_A$ and
$Ann(\gamma_A, \gamma_B)$ both meet $J(G)$, and both $A$ and $B$ do
not meet $J(G)$, there must exists maps $f, g \in G$ such that $J(f)
< \gamma_A$ and $\gamma_A < J(g) < \gamma_B$. Note then that $J(g) <
B$ since $J(g) \cap B =\emptyset$. Since $J(f)$ and $J(g)$ lie
indifferent components of $J(G)$ (separated by $A$), $J(g) \subset
\overline{\cup_{n \in \N} g^{-n}(J(f))}$, and each $ g^{-n}(J(f))
\cap A =\emptyset$, there exists $n_0 \in \N$ such that $
g^{-n_0}(J(f)) > \gamma_A$ and thus $g^{-n_0}(J(f)) > A$.  By
Lemma~\ref{J(k)} there exists $k \in G$ such that $A< g^{-n_0}(J(f))
< J(k) < g^{-(n_0+1)}(J(f)) < J(g)$.

We now find a sub-semigroup $H'$ that satisfies conclusions (1) -
(4) of the theorem.  Keeping Lemma~\ref{attr-basin} in mind, we see
that we may choose $m_1, m_2 \in \N$ large (as in
Example~\ref{cantor-ex}), such that $\beta_1=k^{m_1}$ and
$\beta_2=g^{m_2}$ generate a sub-semigroup $H'$ of $G$ where $J(H')$
is disconnected and contained in $\overline{Ann(J(k),J(g))}$.
Further, $H'$ is hyperbolic since $P^*(H') \subset P^*(G) \subset
K(f) < J(k)$.  By choosing $U$ to be a suitable open set containing
$\overline{Ann(J(k),J(g))}$ we see that $H'$ satisfies parts (1) and
(2) of the theorem.

By Theorem 2.14(2) in~\cite{Su2}, the hyperbolicity of $H'$ implies
$J(H')=\cup_{x \in \Sigma_2'} J_x$, where $\Sigma_2'$ is the
sequence space corresponding to the maps $\beta_1$ and $\beta_2$.
The fact that $J_{x_1} \neq J_{x_2}$ when $x_1 \neq x_2$ follows in
much the same way as the proof that the standard middle-third Cantor
set is totally disconnected.  We present the details now.  First we
define $\sigma$ to be the shift map on $\Sigma_2'$ given by
$\sigma(\gamma_1, \gamma_2, \dots ) = (\gamma_2, \gamma_3, \dots)$.
Then, for $x=(\gamma_1,\gamma_2, \dots )$, one can show by using the
definition of normality $J_x=\gamma_1^{-1}(J_{\sigma(x)})$ and thus
by induction $J_x=\gamma_1^{-1}(J_{\sigma(x)})= \dots = (\gamma_n
\circ \dots \circ \gamma_1)^{-1}(J_{\sigma^n(x)}) \subset(\gamma_n
\circ \dots \circ \gamma_1)^{-1}(J(H'))$.  Thus $J_x \subset
\cap_{n=1}^\infty (\gamma_n \circ \dots \circ
\gamma_1)^{-1}(J(H'))$. But by (induction on) condition (1) we can
see that this intersection will produce distinct sets for distinct
sequences in $\Sigma'$.  Thus we have shown that $J_{x_1} \neq
J_{x_2}$ when $x_1 \neq x_2$

Each $J_x$ is connected by Lemma 3.6 in~\cite{Su03}.  Hence we have
shown parts 3(a) and 3(b).  Now part (4) is then clear by 3(a),
3(b), and Theorem~\ref{mainth1}\eqref{mainth1-1}.  Part 3(c) now
follows directly from Proposition~\ref{qc-surg}.

We have thus shown that $H'=\langle \beta_1, \beta_2 \rangle$
satisfies items (1) - (4) of the theorem, but it is not certain that
$J(g)$ does not meet $\overline{B}$, and so (5) remains in question.
However, letting $\alpha_1=\beta_1$ and $\alpha_2=\beta_1 \circ
\beta_2$ (note $J(\alpha_2)<J(\beta_2) < B$), we see that $H=\langle
\alpha_1, \alpha_2 \rangle$ will satisfy (1) - (5).  We have thus
proved the result in the case that $A$ and $B$ are both doubly
connected Fatou components.

Consider the case where $B$ is the unbounded component of $F(G)$ containing
$\infty$.  As above we obtain $f, g \in G$ such that $J(f)<\gamma_A<
J(g)$ where $\gamma_A$ is a non-trivial curve in $A$.  We then
follow the above method to complete the proof.

Finally, we consider the case where $B=\K$.  As above we obtain
$f, g \in G$ such that $J(g)<\gamma_A< J(f)$ where $\gamma_A$ is a
non-trivial curve in $A$.  We then follow the above method, noting
that the surrounding order inequalities are now reversed from above,
to complete the proof.
\end{proof}

\section{Structural properties of $\J$}\label{struct-prop}
In this section we discuss issues related to the topological nature
of $\J$ as well as discuss issues related to the question of where
the ``small" Julia sets $J(g)$ for $g \in G$ reside inside of the
larger Julia set $J(G)$. In particular, we investigate the question
of when it is the case that a given $J \in \J$ must contain $J(g)$
for some $g \in G$. Since $\Jmin$ and $\Jmax$ play special roles we
will be particularly interested in when these components of $J(G)$
have this property.  When $G= \langle g_\la : \la \in \La \rangle$,
it is of particular interest to know which $J \in \J$ meet
$J(g_\la)$ for some $\la \in \La$.  The first result in this
direction is the following, which appears as Proposition 2.24
in~\cite{Su01}.

\begin{proposition}\label{compact-gen} If $G \in \G$ is generated by a compact family in
Poly, then both $\Jmin$ and $\Jmax$ must contain the Julia set of
one of the generating maps of $G$.
\end{proposition}

In order to succinctly discuss such issues we make the following
definitions.
\begin{definition}
Let $G=\langle h_\lambda: \lambda \in \Lambda \rangle \in
\mathcal{G}$. We say that $J \in \J$ has property $(\star)$ if $J$
contains $J(g)$ for some $g \in G$.  We say that $J \in \J$ has
property $(\star\lambda)$ if $J$ contains $J(h_\lambda)$ for some
generator $h_\lambda \in G$.
\end{definition}

\begin{remark}
A given rational semigroup $G$ may have multiple generating sets.
For example, the whole semigroup itself can always be taken as a
generating set.  However, in this paper when it is written that
$G=\langle h_\lambda:\la \in \La \rangle$, it is assumed that this
generating set is fixed and thus the property $\stl$ is always in
relation to this given generating set.
\end{remark}

\begin{lemma}\label{star-starl}
Let $G=\langle h_\lambda: \lambda \in \Lambda \rangle \in \mathcal{G}$.  \\
a) If $\Jmin$ has property $(\star)$, then $\Jmin$ has property $(\star\lambda)$.\\
b) If $\Jmax$ has property $(\star)$, then $\Jmax$ has property
$(\star\lambda)$.
\end{lemma}

\begin{remark}
Lemma~\ref{star-starl} does not apply to a general $J \in \J$.
Indeed, in Example~\ref{cantor-ex} we see that only $\Jmin$ and
$\Jmax$ have property $\stl$, although infinitely many other $J \in
\J$ have property $(\star )$.
\end{remark}

\begin{proof}
Suppose $J(g) \subset \Jmin$ for $g=h_{\la_1} \circ \dots \circ
h_{\la_k}$ and $\Jmin \cap J(h_{\la}) = \emptyset$ for all $\la \in
\La$.  By Corollary~\ref{comp-surr} we have $\Jmin < J(h_{\la})$ for
all $\lambda \in \Lambda$, and thus by Lemma~\ref{pullback-order}, 
Lemma~\ref{connect-to-Jf} and Theorem~\ref{Jmin} we have $\Jmin < h_{\la}^{-1}(\Jmin)$.
So it also follows from Corollary~\ref{comp-surr} that $J(g) <
h_{\la}^{-1}(J(g))$ for all $\lambda \in \Lambda$. Thus $J(g) <
h_{\la_1}^{-1}(J(g))$ and by Lemma~\ref{backimage}
$J(g)<h_{\la_2}^{-1}(J(g))<h_{\la_2}^{-1}h_{\la_1}^{-1}(J(g))$. By
repeated application of this argument we then get that
$J(g)<h_{\la_k}^{-1}(J(g))<\dots <h_{\la_k}^{-1}\dots
h_{\la_1}^{-1}(J(g))=g^{-1}(J(g))=J(g)$, which is a contradiction.
From this part (a) follows.

Part (b) follows in a similar manner.
\end{proof}

\begin{corollary}
If $\Jmin$ (respectively, $\Jmax$) has non-empty interior, then
$\Jmin$ (respectively, $\Jmax$) has property $\stl$.
\end{corollary}

\begin{proof}
Suppose $int (\Jmin) \neq \emptyset$.  Since $J(G)=\overline{\cup_{g
\in G} J(g)}$ some $J(g)$ must meet $int (\Jmin)$.  Thus it follows
from Lemma~\ref{star-starl} that $\Jmin$ has property $\stl$.
\end{proof}
It is not always the case, however, that  $\Jmin$ and $\Jmax$ have
property $\stl$.

\begin{example}\label{no-lambda-star}
We will give an example of an infinitely generated $G \in \G_{dis}$
such that

\begin{enumerate}
\item $\Jmax$ does not have property $\stl$,

\item $\# \J = \aleph_0$, and

\item there exists $J', J'' \in \J$ and $g \in G$ such that $J'<J''$ and
$g^*(J')=g^*(J'')$.
\end{enumerate}

Set $b_n=2-1/n$ for $n \in \N$ and
$\epsilon_n=\min\{\frac{b_{n+1}-b_n}{10},
\frac{b_{n}-b_{n-1}}{10}\}.$  \\
Set $A=\overline{Ann(0;1/2,2)},
A_n=Ann(0;b_n-\epsilon_n,b_n+\epsilon_n)$ and
$A_n'=Ann(\epsilon_n;b_n-\epsilon_n,b_n+\epsilon_n)$.  Note that the
set $\tilde{A_n}:=\overline{A_n \cup A_n'} \subset
\overline{Ann(0;b_n-2\epsilon_n,b_n+2\epsilon_n)}$ and so, by the
choice of the $\epsilon_n$, the $\tilde{A_n}$ are disjoint. Choose
polynomials $f_n$ and $g_n$ such that $J(f_n)=C(0,b_n)$ and
$J(g_n)=C(\epsilon_n,b_n)$. Choose $m_n \in \N$ large enough so that
$h_n=f_n^{m_n}$ yields $h_n^{-1}(A) \subset A_n,
h_n(\overline{B(0,1/2)}) \subset B(0,1/2)$ and $h_n(\{|z|\geq
2\})\subset \{|z|>5\}$.  Choose $j_n \in \N$ large enough so that
$k_n=g_n^{j_n}$ yields $k_n^{-1}(A) \subset A_n',
k_n(\overline{B(0,1/2)}) \subset B(0,1/2)$ and $k_n(\{|z|\geq
2\})\subset \{|z|>5\}$.  Note that
$\tilde{A}:=\overline{\cup_{n=1}^\infty \tilde{A_n}} \subset A$. Let
$G=\langle h_n, k_n: n \in \N\rangle$ and note that $P^*(G) \subset
B(0,1/2)$ and $G^{-1}(\tilde{A}) \subset G^{-1}(A) \subset
\tilde{A}$ which implies $J(G) \subset\overline{
G^{-1}(\tilde{A})}\subset \tilde{A}$.

We see by forward invariance that $\{|z|> 2\} \subset F(G)$, but
note that $C(0,2) \in \J$ (since the open sets
$Ann(\tilde{A}_n,\tilde{A}_{n+1})$ are all in $F(G)$). Also, for no
$g \in G$ does $J(g)$ meet $C(0,2)$ else there would exist $z_0 \in
C(0,2) \cap J(g)$ such that $|g(z_0)| \leq 2$, contradicting the
fact that each $g \in G$ maps $C(0,2)$ into $\{|z|>5\}$.  Thus
$\Jmax=C(0,2)$ fails to have property $(\star)$.

We now show that $\#\J = \aleph_0$. Letting $J_n \in \J$ be such
that $J_n$ contains the overlapping sets $J(h_n)$ and $J(k_n)$ we
note that, since $J_n \subset \tilde{A}_n$ and the $\tilde{A_n}$ are
separated from each other, each $J_n$ is isolated from the other
$J_m$, i.e., for each $n$ there exists $\epsilon_n
>0$ such that the $\epsilon_n$-neighborhood $B(J_n, \epsilon_n)$
does not meet any other $J_m \in \J$.

Let $C=C(0,2)$.  Since $C \subset J(G)$ we see that
$J(G)=\overline{G^{-1}(C)}$.  We now show that for each $g \in G$,
the set $g^{-1}(C) \subset J_n$ for some $n$.  Write $g=g_{i_1}
\circ \dots g_{i_j}$ where each $g_{i_\ell}$ is a generator for $G$.
Suppose that $g_{i_j}=h_n$ for some $n$.  Then, by the backward
invariance of $A$ under each map in $G$, we have that $g^{-1}(C)
\subset g^{-1}(A) \subset h_n^{-1}(A) \subset A_n$. Since
$g^{-1}(C)$ is connected, $g^{-1}(C)> B(0,1/2)$, and $J(k_n) \cup
J(h_n)$ meets both the inner boundary and outer boundary of $A_n$,
we must have that $g^{-1}(C)$ meets $J(k_n) \cup J(h_n)$ and thus
$g^{-1}(C) \subset J_n.$  Note that the same argument (using $A_n'$)
holds if $g_{i_j}=k_n$. Thus we have shown that $G^{-1}(C) \subset
\cup_{n \in \N} J_n$. Since the $J_n$ are isolated from each other
and accumulate only to $C$, it follows that
$J(G)=\overline{G^{-1}(C)} \subset \overline{\cup_{n \in \N} J_n}
\subset C \cup \bigcup_{n \in \N} J_n$ and so $\#\J = \aleph_0$.
Note also then that we must have $\Jmin = J_1$ and so $\Jmin$ does
have property $\stl$.

We now show that $J'<J''$ for $J', J'' \in \J$ does not necessarily imply
$g^*(J')<g^*(J'')$.  Indeed, we see that $h_n^*(J) = J_n$ for all $J
\in \J$.

We note that we could easily adapt this example (by letting $b_n=
0.5 + 1/n$) to produce $G_1 \in \G$ such that $\Jmin(G_1)$ does not
have property $\stl$, but $\Jmax$ does.  Or we could produce
$G_2=\langle G,G_1 \rangle \in \G$ such that neither $\Jmin$ nor
$\Jmax$ has property $\stl$.
\end{example}

Note that in the above example(s) where $\Jmin$ (respectively
$\Jmax$) did not meet $\cup_{\la \in \La} J(g_\la)$, it was true
that $\Jmin$ (respectively $\Jmax$) was contained in
$\overline{\cup_{\la \in \La} J(g_\la)}$.  We will prove in
Theorem~\ref{almost-lambda-star} that this is indeed always the
case. First we need to prove the following lemma.

\begin{lemma}\label{M'-M''}
Let $\{ g_{\lambda }\} _{\lambda \in \Lambda }\subset G\in \mathcal{G}_{dis}.$
Let $\mathcal{C}$ denote the connected components of
$\overline{\cup_{\la \in \La} J(g_\la)}$.
Then both $M'=\min_{C \in
\mathcal{C}} C$ and $M''=\max_{C \in \mathcal{C}} C$ exist (with
respect to the surrounding order $\leq$).  Also, $PH(C) \supset
\K(G) \supset P^*(G)$ for each $C \in \mathcal{C}$.
\end{lemma}

\begin{proof}
First we note that $\infty \in F(G)$ by Theorem~\ref{Jmin} and so
all sets in $\mathcal{C}$ are contained in $\C$.  Let $C \in
\mathcal{C}$.  Suppose that $z_0 \in \K(G) \setminus PH(C)$.  Let
$\gamma$ be a curve in $\CC \setminus PH(C)$ connecting $z_0$ to
$\infty$ and set $\epsilon = dist(\gamma, PH(C))$.  By
Corollary~\ref{isol}, there exists $\la \in \La$ such that $J(g_\la)
\subset B(C, \epsilon)$.  Hence, we see that $\gamma$ is outside
$J(g_\la)$ implying that $g_\la^n(z_0) \to \infty$ and thus
contradicting the fact that $z_0 \in \K(G)$.

Lemma~\ref{surr-order} shows that the compact connected sets in
$\mathcal{C}$ are linearly ordered with respect to the surrounding
order.  The existence of $M'$ and $M''$ then follows directly from
Lemma~\ref{min-max-exist}.
\end{proof}

\begin{theorem}\label{almost-lambda-star}
Consider $G = \langle g_{\la}: \la \in \La \rangle \in \G_{dis}$.
Let $A= \cup_{\la \in \La} J(g_\la)$ and denote by $M'$ and $M''$
the minimal and maximal connected components of $\overline{A}$,
respectively.  Then both $\Jmin \supset M'$ and $\Jmax \supset M''$
and, in particular, both $\Jmin \cap \overline{A} \neq \emptyset$
and $\Jmax \cap \overline{A} \neq \emptyset$. Furthermore, we have
the following.
\begin{enumerate}
\item If $\Jmin \cap A = \emptyset$ (i.e., $\Jmin$ does not have
property $\stl$), then $\Jmin = M'$ and $\Jmin$ is the boundary of
the unbounded component of $\C \setminus \Jmin$.

\item If $\Jmax \cap A = \emptyset$ (i.e., $\Jmax$ does not have
property $\stl$), then $\Jmax =M''$ and $\Jmax$ is the boundary of
the bounded component of $\C \setminus \Jmax$ which contains
$\Jmin$.
\end{enumerate}
\end{theorem}

\begin{remark}
In the above theorem, if $J(G)$ is connected , then
$\Jmin=J(G)=\Jmax$ meets all $J(g)$ such $g \in G$ and thus meets
$A$.
\end{remark}

\noindent \textbf{Open Question:}  We notice in
Example~\ref{no-lambda-star} that $\Jmax \cap A = \emptyset$ and
$\Jmax$ is a simple closed curve.  However, it is not clear, in
general, whether the hypothesis $\Jmax \cap A = \emptyset$ for $G =
\langle g_{\la}: \la \in \La \rangle \in \G_{dis}$, must
\textit{necessarily} lead to the conclusion that $\Jmax$ is a simple
closed curve.  It is also not clear under this hypothesis whether
$\Jmax$ must be the common boundary of exactly two complementary
domains. So we state these as open questions (noting the
corresponding questions regarding $\Jmin$ are also open).

\begin{remark}
We note that we could also use Theorem~\ref{almost-lambda-star} to
see that $\Jmin$ in Example~\ref{no-lambda-star} has property
$\stl$, since $\Jmin$ must meet $\overline{\cup_{j=1}^\infty J(k_j)
\cup J(h_j)}=C(0,2) \cup \cup_{j=1}^\infty (J(k_j) \cup J(h_j))$,
but does not meet $C(0,2)$.
\end{remark}

\begin{proof}
Let $J', J'' \in \J$ be such that $M' \subset J'$ and $M'' \subset
J''$.  Fix $\la \in \La$.  By the minimality of $M'$, we have either
$J(g_\la) \subset J'$ or $J(g_\la)> J'$. Then $g_\la^*(J') \geq J'$
by Lemmas~\ref{pullback-order} and~\ref{connect-to-Jf} and Theorem~\ref{Jmin}.  Since
$g_\la^*(J') \geq J'$ for all $\la \in \La$ we must have that
$J'=\Jmin$ (because the closed set $\cup_{\{J \in \J: J \geq J'\}}
J$ is then backward invariant under each $g \in G$). Similarly we
see that $J''=\Jmax$.  Thus both $\Jmin \cap \overline{A} \neq
\emptyset$ and $\Jmax \cap \overline{A} \neq \emptyset$.

We now prove (2) by first showing that $J(G) \subset \Jmin \cup M''
\cup Ann(\Jmin, M'')$.  Fix $\la \in \La$.  Since $PH(M'') \supset
P^*(G)$ by Lemma~\ref{M'-M''}, we see that $g_\la^{-1}(M'')$ is
connected using Lemma~\ref{backimage}. Thus $g_\la^{-1}(M'') \cap
M''=\emptyset$, else $J(g_\la)$ meets $\Jmax$ by
Lemma~\ref{connect-to-Jf} which violates our hypothesis that $\Jmax
\cap A = \emptyset$. Hence $g_\la^{-1}(M'') < M''$ by
Lemma~\ref{pullback-order}(b).

From the facts that $g_\la^{-1}(\Jmin) \ge \Jmin$ and
$g_\la^{-1}(M'') < M''$ for all $\la \in \La$, we deduce from
Lemma~\ref{backimage} that the closed annulus-type region $A_1=\Jmin
\cup M'' \cup Ann(\Jmin, M'')$ is backward invariant under each
generator (and thus under each $g \in G$).  Hence $J(G) \subset
\Jmin \cup M'' \cup Ann(\Jmin, M'')$ as desired.

We now suppose there exists $w \in \Jmax \setminus M''$.  Such a
point $w$ must necessarily then lie in $Ann(\Jmin, M'')$ (since $w
\in \Jmin$ would imply $\Jmin=\Jmax=J(G)$ and thus $\Jmax$ clearly
meets $A$). Let $U$ be the connected component of $\C \setminus M''$
which contains $\Jmin$.  Note that $w \in U$ by definition of
$Ann(\Jmin, M'')$.  Recall that $\partial \K(G) \subset \Jmin$.  Let
$\gamma$ be a curve in $U$ which connects $w$ to some $z_0 \in
\K(G)$ and set $\epsilon = dist(\gamma, M'')>0$.  By
Corollary~\ref{isol} there exists a generator $g_{\lambda_0} \in G$
such that $J(g_{\lambda_0}) \subset B(M'',\epsilon)$.  Thus $\gamma
\cap J(g_{\lambda_0}) = \emptyset$.  Since $\K(G) \subset
PH(J(g_{\lambda_0}))$, we see that $z_0 \in \gamma \cap
PH(J(g_{\lambda_0}))$ and so $\gamma < J(g_{\lambda_0})$.  Hence
$\{w\} < J(g_{\lambda_0})$ which implies (by
Corollary~\ref{comp-surr}) either $\Jmax < J(g_{\lambda_0})$ or
$\Jmax \cap J(g_{\lambda_0}) \neq \emptyset$.  Since neither of
these can occur we conclude that no such $w$ exists and thus $\Jmax
= M''$.

Recall that $U$ is the bounded component of $\C \setminus \Jmax$
which contains $\Jmin$.  Since $\Jmax \cap A = \emptyset$, we have
that for every $\lambda \in \Lambda$, the set $J(g_\lambda)$ is
contained in $U$. Hence $\overline{A} \subset \overline{U}$ and so
$\Jmax = M'' \subset \overline{A} \subset \overline{U}$, which
implies $\Jmax =\partial U$.

The proof for case (1) is similar, but simpler.  In this case the
point $\infty$ can play the role of $z_0$ in order to help
demonstrate that any point in $\Jmin \setminus M'$ must lie
``outside" of some $J(g_{\lambda})$ (which is a contradiction). We
omit the details.
\end{proof}

\begin{corollary}
When $G = \langle g_{\la}: \la \in \La \rangle \in \G$ with
$\overline{\cup_{\la \in \La} J(g_\la)}=\cup_{\la \in \La}
J(g_\la)$, then both $\Jmin$ and $\Jmax$ must have property $\stl$.
\end{corollary}

\begin{remark}
The above corollary applies, for example, when $G = \langle g_{\la}:
\la \in \La \rangle \in \G_{dis}$ has $\cup_{\la \in \La}
J(g_\la)=\cup_{k=1\dots n} J(g_{\la_k})$.  Such non compactly
generated examples can easily be constructed.  Other more ``exotic"
examples can also be constructed to satisfy the hypotheses of the
corollary.
\end{remark}

\begin{example}
We note that without the hypothesis that $\Jmin \cap A =\emptyset$
in Theorem~\ref{almost-lambda-star}(1), the conclusion that $\Jmin
=M'$ might not hold.  Set $f_1(z)=z^2,
f_2(z)=(z-\epsilon)^2+\epsilon$, and $f_3(z)=z^2/4$. For $\epsilon
>0$ small and $m_1, m_2, m_3$ all large we set
$g_1=f_1^{m_1}, g_2=f_2^{m_2},$ and $g_3=f_3^{m_3}$ and note that
$G=\langle g_1, g_2, g_3 \rangle \in \G_{dis}$.  Thus $M'=J(f_1)
\cup J(f_2)= C(0,1) \cup C(\epsilon, 1)$. However, the real point in
$g_1^{-1}(\{1+\epsilon\})$ is clearly in $g_1^{-1}(J(f_2)) \subset
\Jmin$, but not in $M'$.
\end{example}

\begin{theorem}
Let $G = \langle g_{\la}: \la \in \La \rangle \in \G_{dis}$ and
suppose $\Jmin \cap \cup_{\la \in \La} J(g_\la)= \emptyset$.  Then
$\partial \K(G) \subset M' \subset \overline{\cup_{\la \in \La}
J(g_\la)},$ where $M'$ is the minimal connected component of
$\overline{\cup_{\la \in \La} J(g_\la)}$
\end{theorem}

\begin{proof}
By Theorem~\ref{almost-lambda-star}(1), we see that $\partial \K(G)
\subset J(G) \subset M' \cup \Jmax \cup Ann(M',\Jmax)$.  Moreover,
by Lemma~\ref{M'-M''}, we have $\partial \K(G) \subset PH(M')$.
Hence $\partial \K(G) \subset M'.$
\end{proof}

Having discussed properties $(\star)$ and $\stl$ with respect to
$\Jmin$ and $\Jmax$ we now turn our attention to a general $J \in
\J$.  In particular, we investigate what can be said about which $J
\in \J$ have property $(\star)$ or $\stl$.  We also concern
ourselves with the question of when does every $J \in \J$ have
property $(\star)$ or $\stl$. We begin with the following
definition.

\begin{definition}
Let $G \in \G$.  We say that $J \in \J$ is \textbf{isolated} in $\J$
if there exists $\epsilon >0$ such that $B(J, \epsilon)$ does not
meet any other set in $\J$.
\end{definition}

\begin{lemma}\label{J-isol-star}
Let $G \in \G$ with $J \in \J$ isolated in $\J$.  Then $J$ has
property $(\star)$.
\end{lemma}

\begin{proof}
Assume that $\epsilon >0$ is such that $B(J, \epsilon)$ does not
meet any other set in $\J$.  Since $J(G)=\overline{\cup_{g \in
G}J(g)}$ by Theorem~\ref{HMJ}, we see that any point in $J$ must
have, within a distance $\epsilon$, a point in some $J(g)$, where $g
\in G$.  It must then be the case that $J(g)$, which lies in some
set in $\J$, must lie entirely in $J$.
\end{proof}

\begin{remark}
If $G \in \G$ is such that $\#\J<+\infty$, then clearly each $J \in
\J$ is isolated in $\J$ and so each $J \in \J$ has property
$(\star)$.  We note, however, that if each $J \in \J$ is isolated in
$\J$, then it is not necessarily the case that each $J \in \J$ has
property $\stl$.  See the proof of Theorem~\ref{num-J-comp=k} where,
for any positive integer $k$, a semigroup $G' \in \G$ is constructed
such that $\#\J=k$, but only $\Jmin$ and $\Jmax$ have properly
$\stl$.
\end{remark}

\begin{remark}
If $G=\langle h_\lambda: \lambda \in \Lambda \rangle \in
\mathcal{G}$ where $\#\Lambda \leq \aleph_0$ and $\#\J$ is
uncountable, then since $\#G=\aleph_0$ we see that some $J \in \J$
must fail to have property $(\star)$.  An example of this is the
Cantor set of circles in Example~\ref{cantor-ex}.
\end{remark}

\begin{example}
Suppose $G \in \mathcal{G}$ and $\#\J = \aleph_0$.  Then it is
possible that not all $J \in \J$ have property $(\star)$ as in
Example~\ref{no-lambda-star}.  But it is also possible that all $J
\in \J$ do have property $(\star)$ as in~\cite{Su01} Theorem 2.26
where both $\Jmin$ and $\Jmax$ have property $\stl$ and each other
component of $J(G)$ is isolated in $\J$.
\end{example}

We saw above that isolated $J \in \J$ have property $(\star)$.  We
now show that this is also true for the components of $J(G)$ which
contain the pre-image of an isolated $J \in \J$.

\begin{claim}\label{back-isol-star}
Let $G \in \G$.  If $J_1 \in \mathcal{J}$ is isolated in $\J$ and
$h^{-1}(J_1) \subset J \in \mathcal{J}$ for some $h \in G$, then $J$
has property $(\star)$.
\end{claim}

\begin{proof}
Since $J_1$ is isolated in $\J$, Lemma~\ref{J-isol-star} implies
there exists $g \in G$ such that $J(g) \subset J_1$.  Thus, since
$J_1$ is isolated in $\J$, for large $n \in \N$ we have $g^{-n}(J)
\subset J_1$. Hence $h^{-1}(g^{-n}(J)) \subset h^{-1}(J_1) \subset
J$ which implies $J(g^n\circ h) \subset J$.
\end{proof}

\noindent \textbf{Open Question:}  If $\#\J = \aleph_0$ and $G \in
\G$ is finitely generated, then must every $J \in \mathcal{J}$  have
property $(\star)$? Note that the finitely generated condition is
required by Example~\ref{no-lambda-star}.  Also, Example~\ref{cantor-ex}
shows that if $\#\J$ is uncountable, $\J$ can have (uncountably many) $J$
which fail to have property $(\star)$.\\

We now turn our attention to considering those semigroups where
$\Jmin$ has property $\stl$.  In particular, we examine the
generating maps whose Julia sets meet $\Jmin$ as well as the
sub-semigroup generated by just these special maps.

\begin{definition}
Let $G=\langle h_\lambda: \lambda \in \Lambda \rangle \in
{\G}_{dis}.$  We set
$$ B_{\min }=B_{\min}(G):=\{ \la \in \La : J(h_{\la})\subset J_{\min }(G)\} .$$
If $B_{\min }\neq \emptyset $, let $H_{\min }(G)$ be the sub-semigroup of $G$ which is
generated by $\{ h_{\la}: \la \in B_{\min }\} .$
\end{definition}

The following proposition has been shown in \cite{Su01, Su9}.

\begin{proposition}
\label{bminprop}
Let $G=\langle h_\lambda: \lambda \in \Lambda \rangle \in
{\G}_{dis}.$ If $\{h_\lambda:\lambda \in \Lambda\}$ is compact in
$\textrm{Poly}$, then $B_{\min }$ is a proper non-empty subset of
$\La$ under the above notation.
\end{proposition}

\begin{proof}
The result follows since by Proposition~\ref{compact-gen} both
$\Jmin$ and $\Jmax$ have property $\stl$.
\end{proof}

It is natural to investigate the relationship between $H_{\min}(G)$
and $G$.  Specifically we ask, and answer, the following questions
for a semigroup $G \in \G_{dis}$:
\vskip.1in
\begin{tabular}{lll}
  &(1) Must $J(\Hmin(G))=\Jmin(G)$? & (2) Must $J(\Hmin(G))$ be
connected? \\
  &(3) Must $\Jmin(H_{\min}(G))=\Jmin(G)$? & (4) Must
$\Hmin(\Hmin(G))=\Hmin(G)$?
\end{tabular}
\vskip.1in

\noindent The answer to each of these questions is NO, as we see in
this next example.

\begin{example}\label{Hmin-not}
We will construct a single 3-generator polynomial semigroup $G \in
\G_{dis}$ which negatively answers questions (1)-(4). Furthermore,
we will show that $\# \J = \aleph_0$.

Let $h_1(z)=-z^2$ and $f_2(z)=z^2/\sqrt{2}$ and note that
$J(h_1)=C(0,1)$ and $J(f_2)=C(0,\sqrt{2})$.  We set $h_2=f_2^{m_2}$
where conditions on the large $m_2 \in \N$ will be specified later.
Choose point $P \in h_2^{-1}(J(h_1))$ such that $P>0$. Note that
$P=\sqrt{2}-\delta$ where $\delta$ is small for $m_2$ large. Hence
$P^2=2-2\sqrt{2}\delta+\delta^2$ and
$P^4=4-8\sqrt{2}\delta+O(\delta^2)$.

Setting $r=\frac{P^4+P^2}{2}=3-M_1\delta+O(\delta^2)$ and
$\epsilon=\frac{P^4-P^2}{2}=1-M_2\delta+O(\delta^2)$ where $M_1,
M_2>0$, we see that for $f_3(z)=\frac{(z-\epsilon)^2}{r}+\epsilon$
both $h_1(P)=-P^2=\epsilon-r$ and $h_1^2(e^{i\pi/4}
P)=P^4=r+\epsilon$ lie in $C(\epsilon,r)=J(f_3)$.

Suppose there exists $w \in h_1^{-1}(J(f_3))\cap J(f_3)$, i.e.,
$h_1(w) \in J(f_3)$ and $w \in J(f_3)$.  Then $|h_1(w)|=|w|^2 \geq
|-P^2|^2=P^4$ since $-P^2$ is the point in $J(f_3)$ of smallest
modulus.  Since $P^4$ is the point in $J(f_3)$ of largest modulus,
we see that $h_1(w)$ could only be in $J(f_3)$ if $w=-P^2$.  But
$h_1(-P^2)=-P^4 \notin J(f_3)$, and so we have $h_1^{-1}(J(f_3))\cap
J(f_3) = \emptyset$ and thus $h_1^{-1}(J(f_3))<J(f_3)$ (see
Figure~1).

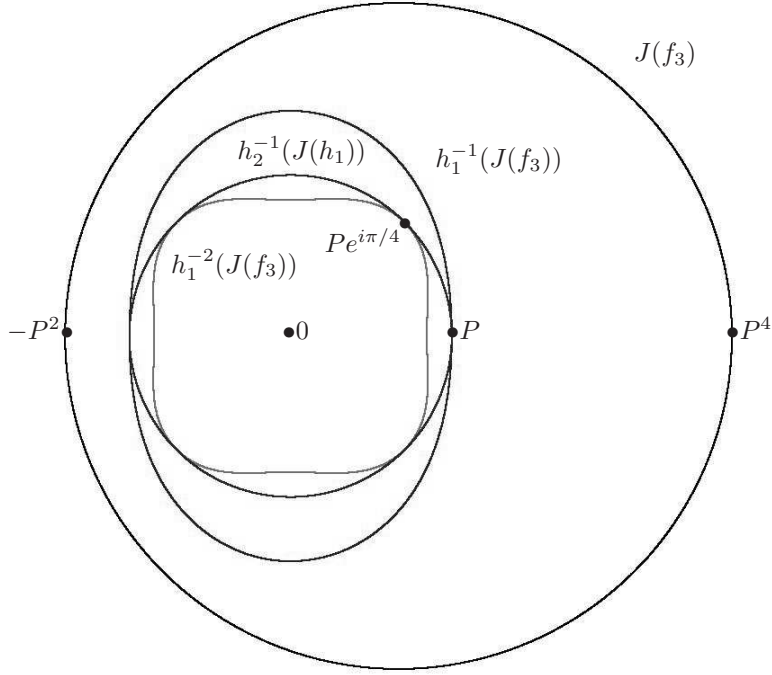
\begin{figure}
\begin{center}
\epsfig{file=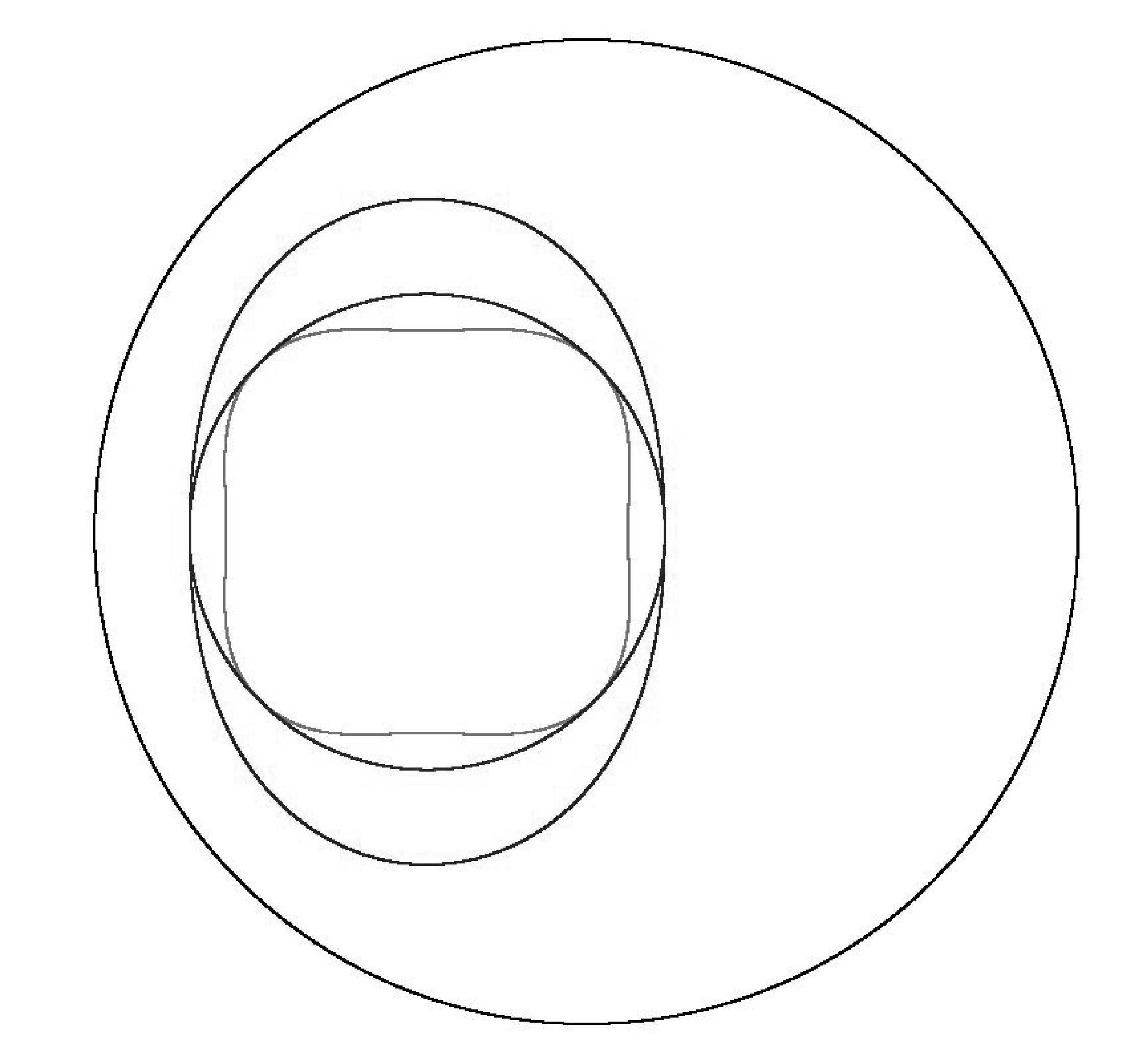,width=4in}
\put(-125,200){$h_1^{-1}(J(f_3))$} \put(-50,240){$J(f_3)$}
\put(-225,160){$h_1^{-2}(J(f_3))$}
\put(-200,203){$h_2^{-1}(J(h_1))$} \put(-121,135){$\bullet P$}
\put(-15,135){$\bullet P^4$} \put(-183,135){$\bullet 0$}
\put(-287,135){$-P^2 \bullet$} \put(-167,168){$ Pe^{i \pi/4}$}
\put(-139,176){$ \bullet$} \caption{This figure shows the sets (from
outside to inside) $J(f_3), h_1^{-1}(J(f_3)), h_2^{-1}(J(h_1)),$ and
$h_1^{-2}(J(f_3))$. For this picture we chose $m_2=5$ and $P =
2^{31/64} \approx 1.39898$. }
\end{center}
\end{figure}

Since for $m_2$ large we clearly have $h_2^{-1}(J(f_3))<J(f_3)$, we
are then free to choose $A'$ to be any closed annulus such that $A'
\subset B(\epsilon, r), int(A') \neq \emptyset,$ and
$A'>h_1^{-1}(J(f_3))\cup h_2^{-1}(J(f_3))$.

Choose $m_3 \in \N$ large enough so that $h_3=f_3^{m_3}$ maps $A'$
into $B(0,1)$ (which can be done since the super attracting fixed
point of $f_3$ is $\epsilon \in B(0,1)$).  This then implies that
$B(0,1)$ is forward invariant under each map $h_1, h_2$ and $h_3$
and $P^*(G) \subset B(0,1)=int\hat{K}(G)\subset F(G)$ where
$G=\langle h_1, h_2, h_3 \rangle$.

Let $A=J(h_2)\cup h_2^{-1}(J(h_1))\cup h_1^{-1}(J(h_3)) \subset
J(G)$.  Because $h_1^{-1}(J(h_3))$ is connected
(Corollary~\ref{backcor}) and contains both the point $iP^2$, which
lies outside the circle $J(h_2)$, and the point $P \in
h_2^{-1}(J(h_1))$, which lies inside the circle $J(h_2)$, we see
that the set $A$ is connected.

Note that $e^{i\pi/4} P \in h_2^{-1}(J(h_1)) \subset A$ and
$e^{i\pi/4} P \in h_1^{-1}(h_1^{-1}(J(h_3))) \subset  h_1^{-1}(A)$
since $h_1^2(e^{i\pi/4}P)=P^4 \in J(h_3)$.  Thus $h_1^{-1}(A) \cap A
\neq \emptyset$ and since $h_1^{-n}(A)$ is connected for each $n \in
\N$ by Lemma~\ref{backimage}, we see that Lemma~\ref{connect-to-Jf}
implies $A$ and $J(h_1)$ are contained in a single $J \in \J$. Since
$P^*(G) \subset B(0,1) \subset F(G)$, we see that $J=\Jmin$. Thus
both $J(h_1)$ and $J(h_2)$ are contained in $\Jmin(G)$.  Note that
$\partial \K(G) = J(h_1) \neq \Jmin(G)$.

Both $h_1$ and $h_2$ map $A'$ into the unbounded component of $\CC
\setminus J(h_3)$ (since $A'$ is outside of both $h_1^{-1}(J(h_3))$
and $h_2^{-1}(J(h_3))$), which is forward invariant under each map
$h_1, h_2$ and $h_3$.  The map $h_3$ maps $A'$ into $B(0,1)$, which
is also forward invariant under each map $h_1, h_2$ and $h_3$. Hence
for any $g \in G$ we have that $g(A')\cap A' = \emptyset$ and so
$int(A') \subset F(G)$.  We conclude that $J(h_3)$ is not contained
in $\Jmin$.

Thus we have that  $\Hmin(G)=\langle h_1, h_2\rangle$.  One can
easily show that $J(\Hmin(G))$ is disconnected (Cantor set of
circles) and thus $J(\Hmin(G)) \neq \Jmin (G)$.  Also
$\Hmin(\Hmin(G))=\langle h_1\rangle \neq \langle h_1, h_2
\rangle=\Hmin(G)$ and $\Jmin(G) \neq \Jmin(\Hmin(G))$.

We now show $\#\J_G=\aleph_0$.  Consider the set $B=\Jmin(G) \cup
\bigcup_{n=1}^\infty h_3^{-n}(\Jmin(G))\cup J(h_3) \subset J(G)$,
which clearly contains more than three points.
We have $J(h_{3})\subset J_{\max}(G)$ and $h_{1}^{-1}(J(h_{3}))\subset J_{\min }(G).$
By Theorem~\ref{mainth1} or Corollary~\ref{backcor},
we obtain $h_{1}^{-1}(B)\subset h_{1}^{-1}(J(G))\subset J_{\min }(G)\subset B.$
Similarly, taking $m_{2}$ so large,
we may assume $h_{2}^{-1}(J(h_{3}))\subset J_{\min }(G)$, and it implies
$h_{2}^{-1}(B)\subset h_{2}^{-1}(J(G))\subset J_{\min }(G)\subset B.$
Since $B$ is closed
and backward invariant under each generator of $G$ (and hence under
every $g \in G$), we must have that $B=J(G)$.  Also, since
$h_3^{-1}(\Jmin(G))$ is connected (by Corollary~\ref{backcor}) and
does not meet $\Jmin(G)$, we see that $\Jmin(G) <
h_3^{-1}(\Jmin(G))$. Repeated application of Lemma~\ref{backimage}
shows us that $\Jmin(G) < h_3^{-1}(\Jmin(G)) < h_3^{-2}(\Jmin(G)) <
\dots < h_3^{-n}(\Jmin(G)) < \dots$. From this we may conclude that
$\J=\{\Jmin(G), J(h_3),h_3^{-n}(\Jmin(G)): n \in \N \}$, thus
demonstrating that $\#\J_G=\aleph_0$.
\end{example}

\begin{remark}
Note that Example~\ref{Hmin-not} does not settle (in the negative)
the open question stated above since Claim~\ref{back-isol-star} with
$J_1=\Jmin$ shows that each $J \in \J_{G}$ contains $J(g)$ for some
$g \in G$. One could also note that every $J \in \J_{G}$ other than
$J(h_3)$ is isolated in $\J_{G}$ and so from Lemma~\ref{J-isol-star}
each such $J$ has property $(\star)$.
\end{remark}

Question: Does there exists an example of some $G \in \G$ which can
negatively answer questions (1)-(4) addressed by
Example~\ref{Hmin-not}, but where $\#\J_{G}$ is finite? The answer,
as we see in the next example, is YES.  We will also see that this
example will settle two other questions that naturally arise when
considering the two following results.  In~\cite{Su01, Su9} it was shown
that, for each positive integer $k$, there exists a semigroup $G \in
\G_{dis}$ with $2k$ generators such that $J(G)$ has exactly $k$
components. Furthermore, in~\cite{Su04} it was shown that any
semigroup in $\G$ generated by exactly three elements will have a
Julia set with either one or infinitely many components. Hence we
have the following questions.

(5) What is the fewest number of generators that can produce a
semigroup $G \in \G_{dis}$ with $\# \J =3$?

(6) For fixed integer $k > 3$, what is the fewest number of
generators that can produce a
semigroup $G \in \G_{dis}$ with $\# \J =k$? \\
The answer to both of these questions is four as stated in
Theorem~\ref{num-J-comp=k} whose proof is given now.

\begin{proof}[Proof of Theorem~\ref{num-J-comp=k}]
Fix $k \ge 2$ since the $k=1$ case is trivial.  Let maps $h_1, f_2$,
and $f_3$ and integer $m_2 \in \N$ be defined as in
Example~\ref{Hmin-not}. Again, we set $h_2 = f_2^{m_2}$ and $h_3 =
f_3^{m_3}$ where large $m_3 \in \N$ will be specified to fit the
stipulations given below. Letting $\gamma_1$ denote $C(0,1)$ and
$\gamma_2$ denote the boundary of the unbounded component of $\C
\setminus (h_1^{-1}(J(f_3)) \cup h_2^{-1}(J(f_3)))$, we set
$B=\overline{Ann(\gamma_1, \gamma_2)}$. Let $A'$ be any closed
annulus such that $A' \subset B(\epsilon, r), int(A') \neq
\emptyset,$ and $A'>B$.  We choose $m_3 \in \N$ large enough so that
$h_3^{-1}(B) >A' > B$. Set $G=\langle h_1, h_2, h_3 \rangle$ and
note that for $m_3$ large enough, $\gamma_1 \cup \gamma_2 \subset
\Jmin(G) \subset B$, as in Example~\ref{Hmin-not}.

Set $\ell = k-2$.  Let $C=C(\epsilon_0,r_0)$ be the circle which is
internally tangent to the circle $J(h_3)$ at the point
$\epsilon+r=P^4$ such that $C$ meets $h_3^{-(\ell +1)}(B)$ and
$\overline {B(\epsilon_0, r_0)} \supset h_3^{-(\ell +1)}(B)$. Hence
$C$ must necessarily meet $h_3^{-(\ell +1)}(\gamma_2) \subset
h_3^{-(\ell +1)}(\Jmin(G))$ and $C > h_3^{-\ell}(B) > \dots >
h_3^{-1}(B)>B$. Note that as $m_3 \to \infty$, we have $\epsilon_0
\to \epsilon$ and $r_0 \to r$.  We may assume then that $m_3$ has
been chosen large enough so that $\epsilon_0 \in B(0,1)$.

Set $f_4(z)= \frac{(z-\epsilon_0)^2}{r_0} +\epsilon_0$ and observe
that $J(f_4)=C$.  Let $R>0$ be large enough so that $f_j(\C
\setminus \overline{B(0,R)}) \subset \C \setminus \overline{B(0,R)}$
for $j=1, \dots, 4$ (where $f_1:=h_1$).  Let $A'' $ be a closed
annulus such that $int(A'') \neq \emptyset$ and
$h_{3}^{-\ell}(B)<A''<h_{3}^{-(\ell+1)}(B).$ Then
$h_{3}^{\ell+1}(A'')\subset B(0,1).$ Let $A_{0}:= A'\cup \cup
_{j=0}^{\ell}h_{3}^{j}(A'').$  We define $h_4=f_4^{m_4}$ where $m_4
\in \N$ is large enough such that (i) $h_4(A_0) \subset B(0,1)$,
(ii) $h_4^{-1}(\gamma_1)$ meets $h_3^{-(\ell +1)}(\Jmin(G))$ (this
is possible since the connected set $h_3^{-(\ell +1)}(\Jmin(G))$
meets, but is not contained in, $C$),
(iii) $h_{4}(\cup _{j=1}^{2} \cup _{h\in \langle
h_{1},h_{2},h_{3}\rangle \cup \{ id\} } h\circ h_{j}(A_0))\subset \C
\setminus \overline{B(0,R)}$ (note that $\overline{h_1(A_0) \cup
h_2(A_0)} \subset \CC \setminus K(h_{3})$, which is equal to the
connected component of $F(\langle h_{1},h_{2},h_{3}\rangle )$
containing $\infty $), and (iv) $h_{4}^{-1}(\gamma
_{1})>h_{3}^{-\ell}(B).$

Set $G' = \langle h_1, h_2, h_3, h_4 \rangle$.  Since $B(0,1)$ is
forward invariant under each map in $G'$, we conclude $B(0,1)
\subset F(G')$ and $P^*(G') \subset B(0,1)$. Thus $G' \in \G.$
 Also, $int(A_0) \subset F(G')$ since one can show that
$g(A_{0})\subset A_{0}\cup B(0,1)\cup \hat{\Bbb{C}}\setminus
K(h_{3})$ for all $g \in G'$. Hence $G'\in {\mathcal G}_{dis}.$  By
applying Lemma~\ref{connect-to-Jf} and Lemma~\ref{pullback-order}
(noting that $\Jmin(G') < int(A') < J(h_3)$), we have that
$h_3^{-n}(\Jmin(G')) < h_3^{-(n+1)}(\Jmin(G'))$ for all $n \geq 0$.
Further, $\Jmax(G')$, which must contain $J(h_3)$ and $J(h_4)$ by
Proposition~\ref{compact-gen}, must also contain
$h_3^{-n}(\Jmin(G'))$ for all $n> \ell$. By examining the dynamics
one can then show that $J(G') = \Jmin(G') \cup \Jmax(G') \cup
h_3^{-1}(\Jmin(G')) \cup \dots \cup h_3^{-\ell}(\Jmin(G'))$, since
this set is closed and backward invariant under each generator of
$G'$.  Moreover, since $h_{3}^{-\ell}(\gamma
_{1})<A''<h_{3}^{-(\ell+1)}(\gamma _{1})$, $h_{3}^{-\ell}(\gamma
_{1})\subset h_{3}^{-\ell}(J_{\min }(G'))$,
$h_{3}^{-(\ell+1)}(\gamma _{1})\subset J_{\max}(G')$ and
$int(A'')\subset F(G')$, we have $h_{3}^{-\ell}(J_{\min
}(G'))<int(A'')<J_{\max }(G').$ Thus we see that $J(G')$ has exactly $k$
components.
\end{proof}

\begin{remark}
The addition of one generating function in the proof of
Theorem~\ref{num-J-comp=k} to the semigroup in
Example~\ref{Hmin-not} illustrates something of a general principle
(which we decline to attempt to make precise) at work when dealing
with the dynamics of semigroups in $\G_{dis}$.  Namely, if one adds
a generator (or a whole family of generators) whose Julia set does
not meet $\Jmin$ of the new larger semigroup, then key properties of
the dynamics can often be preserved.  See for example
Theorems~\ref{semi-hyp} and~\ref{hyp}.

However, as we see in this next lemma, adding ``too many" new
functions will necessarily destroy certain critical aspects of the
dynamics. In particular, if we look to produce a new semigroup in
$\G$ by adding ``too many" generating polynomials of small degree
(such that $G_k$ defined in the lemma is not pre-compact) to a
semigroup $G \in \G_{dis}$, then the new semigroup will necessarily
have a \textbf{connected} Julia set.
\end{remark}

\begin{lemma}\label{gen-pre-cpct}
Let $G=\langle h_\lambda: \lambda \in \Lambda \rangle \in
\mathcal{G}_{dis}$.  Then each $G_k=\{g \in G: \deg (g) \leq k\}$ is
pre-compact in {\em Poly} and, in particular, each $\{h_\lambda:\lambda
\in \Lambda\} \cap G_k$ is pre-compact in {\em Poly}.
\end{lemma}

\begin{remark}
As stated earlier, a possibly generating set for $G$ is $G$ itself,
which is necessarily not pre-compact (since it contains elements of
arbitrarily high degree).  Thus it is impossible to strengthen
Lemma~\ref{gen-pre-cpct} to conclude that $\{h_\lambda:\lambda \in
\Lambda\}$ is pre-compact.
\end{remark}

\begin{proof}
Note that $J(G)$ is bounded in $\C$ since Theorem~\ref{Jmin} yields
$\infty \in F(G)$. Choose $R>0$ such that $J(G) \subset
\overline{B(0,R)}$.  Then $\textrm{Cap}(J(g)) \leq R$ for all $g \in
G$, where $\textrm{Cap}(E)$ denotes the logarithmic capacity of the
set $E$ (see~\cite{A} for definition and properties). Also, since $G
\in \G_{dis}$ we have $int(\K(G)) \neq \emptyset$ (see Theorem~\ref{Jmin} or \cite{Su01}),
and so there exists a ball of some radius $r>0$ in $\K(G)$. Thus
$\textrm{Cap}(J(g)) \geq r$ for all $g \in G$.

Let $H_n=\{g \in G: \deg (g) =n\}$. In order to show that
$G_k=\cup_{n=1}^k H_n$ is pre-compact, it suffices to show that each
$H_n$ is pre-compact. We now fix $g(z)=a_n z^n +\dots + a_0$ in
$H_n$ and proceed to show that $|a_n|$ is uniformly bounded below by
$R^{1-n}$ and uniformly bounded above by $r^{1-n}$, and that the
remaining coefficients $a_{n-1}, \dots, a_0$ of $g(z)$ are uniformly
bounded (above) by positive constants which only depend on $r, R$
and $n$. Recalling Remark~\ref{Poly}, it follows then that $H_n$ is
pre-compact. Since $|a_n|^{-1/(n-1)}=\textrm{Cap}(J(g))$
(see~\cite{CG}, p.~35), we see that $r^{1-n} \geq |a_n| \geq
R^{1-n}$.  Express $g'(z)=\beta(z-\alpha_1)\dots(z-\alpha_{n-1})$
where $\beta = na_n$ and the $\alpha_j$ are the critical points of
$g$ which, since $G\in {\mathcal G}$ and $\infty \in F(G)$, must lie
in $\Bbb{C}\setminus F_{\infty }(G)\subset \overline{B(0,R)}$, where
$F_{\infty }(G)$ denotes the connected component of $F(G)$
containing $\infty .$

One can multiply out the terms in the expansion of $g'(z)$ and find
an anti-derivative to see that the $a_{n-1}, \dots, a_1$
coefficients of $g(z)$ are also bounded by constants which depend
only on $r, R$ and $n$. Now fix $z_0 \in \K(G)$. Since
$g(\hat{K}(G)) \subset \hat{K}(G) \subset \overline{B(0,R)}$, we
have $|g(z_0)| = |a_n z_0^n + \dots +a_0| \leq R$.  Thus, since
$|a_{n}|, \dots, |a_1|$ are bounded by constants depending only on
$r, R$ and $n$, the same is true for $|a_0|$.
\end{proof}

\begin{remark}
The proof of Lemma~\ref{gen-pre-cpct} also holds for any $G \in \G$
such that there exists both lower and upper bounds on Cap$(J(g))$
for all $g \in G$ (e.g., when $\K(G)$ contains some non-degenerate
continuum and $\infty \in F(G)$).
\end{remark}

\section{Proof of Theorems~\ref{semi-hyp} and~\ref{hyp}}\label{pf-hyp}

\begin{example}\label{nothyp}
Let $f_1(z)=z^2+c$ where $c>0$ is small (thus $J(f_1)$ is a
quasi-circle).  Let $z_0 \in \R$ denote the finite attracting fixed
point of $f_1$.  Note that $f_1^k(0)$ increases to $z_0$. Choose
$f_2(z)=\frac{(z-z_0)^2}{(c-z_0)}+z_0$ and note that $J(f_2)=C(z_0,
|c-z_0|)$. For $m_1, m_2 \in \N$ large $h_1=f_1^{m_1}$ and
$h_2=f_2^{m_2}$ each map $B(z_0, |c-z_0|)$ into itself and $J(G)$ is
disconnected for $G=\langle h_1, h_2 \rangle$. Note that $P^*(G)
\subset \overline{B(z_0, |c-z_0|)}$ and so $G \in \mathcal{G}$. We
have $H=\langle h_2 \rangle$ is hyperbolic, but since $f_1(0)=c\in
J(h_2) \subset J(G)$, the semigroup $G=\langle H,h_1 \rangle$ is not
hyperbolic even though $J(h_1) \cap \Jmin(G) = \emptyset$.

By conjugating $h_2$ by a suitable rotation we may assume that
$\{h_2^k(c):k \in \N\}$ is dense in $J(h_2)$ and therefore we see
that $H$ can be hyperbolic and have $G$ fail to even be
sub-hyperbolic. However, Theorem~\ref{semi-hyp} does imply that
$G=\langle H,h_1 \rangle$ is semi-hyperbolic.
\end{example}

\begin{remark}\label{ghyp-not-Ghyp}
In contrast to the analogous behavior of Iterated Function Systems
where contraction in each generating map leads to a semigroup (IFS)
that is overall contracting, we see that in Example~\ref{nothyp} each
map of the semigroup $G$ is hyperbolic, yet the entire semigroup $G$
fails to be hyperbolic.  To see this, note that each map $h_2^n$ is
hyperbolic and for each map $g \in G \setminus \{h_2^n\}$ we have
$P^*(g) \subset P^*(G) \subset \overline{B(z_0,|c-z_0|)}$ and $J(g)
>\overline{B(z_0,|c-z_0|)}$ which implies $g$ is hyperbolic.
\end{remark}

We now state a lemma which we will use the proof of
Theorem~\ref{semi-hyp}.

\begin{lemma}~\label{int=int}
Let $H_{1}$ be a polynomial semigroup in $\G$ and let $\Gamma$ be a
compact family in {\em Poly}.  Let $H_{2}=\langle H_1, \Gamma \rangle$ be
the semigroup generated by $H_{1}$ and $\Gamma.$
Suppose \\
\indent (1) $H_{2}\in \G _{dis}$, and \\
\indent (2) $J(\gamma)\cap J_{\min }(H_{2})=\emptyset$ for $\gamma \in \Gamma$. \\
Then $\hat{K}(H_{1})=\hat{K}(H_{2})$, which then implies
$\Jmin(H_1)\subset \Jmin(H_2)$ since $\partial \hat{K}(H_{1})
\subset \Jmin(H_1)$ and $\partial \hat{K}(H_{2}) \subset
\Jmin(H_2)$.
\end{lemma}

\begin{remark}\label{facts}
We recall the facts given in~\cite{Su01, Su9} that for any $G \in \G $ we
have $int \K(G)=\K(G) \cap F(G).$ Moreover, for any $G\in \G_{dis}$,
we have $int \K(G) \neq \emptyset$ and $g(\K(G)\cup J_{\min }(G)) \subset int \K(G)$
for any $g \in G$ such that $J(g) \cap \Jmin(G) = \emptyset$.
\end{remark}

\begin{proof}
We begin by first showing that $J_{\min }(H_{1})\subset J_{\min }(H_{2}).$
Let $\mathcal{C}$ be the set of all connected components of
$\overline{\bigcup _{\gamma \in \Gamma }J(\gamma )}.$
By Lemma~\ref{M'-M''}, $M_{1}:=\min _{C\in \mathcal{C}}C$ exists with respect to
the surrounding order. Let $J_{1}\in \mathcal{J}_{H_{2}}$ be the element containing
$M_{1}.$  Let $J_{2}\in \mathcal{J}_{H_{2}}$ be the element containing $J_{\min }(H_{1}).$
Let $J_{0}:=\min \{ J_{1},J_{2}\} \in \mathcal{J}_{H_{2}}.$
Then for each $g\in H_{1}\cup \Gamma $,
either $J_{0}<J(g)$ or $J(g)\subset J_{0}.$ By Lemma~\ref{pullback-order} and Lemma~\ref{connect-to-Jf},
we obtain that for each $g\in H_{1}\cup \Gamma $, either $g^{-1}(J_{0})>J_{0}$ or
$g^{-1}(J_{0})\subset J_{0}.$ By Corollary~\ref{backcor} or Theorem~\ref{mainth1},
it follows that $A:=\bigcup _{J\in \mathcal{J}_{H_{2}}, J\geq J_{0}}J $
is closed, $\sharp A\geq 3$, and $g^{-1}(A)\subset A$ for each $g\in H_{2}.$
Therefore $J(H_{2})\subset A.$ Thus $J(H_{2})= A$ and hence $J_{0}=J_{\min }(H_{2}).$
From assumption (2), however, it must be the case that $J_{0}=J_{2}.$
Therefore $J_{\min }(H_{1})\subset J_{2}=J_{\min }(H_{2})$ as desired.

By Remark~\ref{facts} and Theorem~\ref{Jmin}, it follows that
for each $\gamma \in \Gamma $, we have
$\gamma (\partial \hat{K}(H_{1}))\subset \gamma (J_{\min }(H_{1})) \subset \gamma (J_{\min } (H_2))
\subset int\hat{K}(H_{2}).$
Let $\gamma \in \Gamma $ and
let $U$ be any connected component of $int\hat{K}(H_{1}),$ which we note is simply connected by the maximum principle.
Then $\partial \gamma (U)\subset \gamma (\partial U) \subset int\hat{K}(H_{2}).$
Let $V$ be the connected component of $int\hat{K}(H_{2})$ containing the connected set $\gamma (\partial U).$
Then $\partial \gamma (U)\subset V.$
By the maximum principle, $\overline{\Bbb{C}}\setminus \overline{V}$ is connected and unbounded.
Hence $\gamma (U)\subset V.$
From this argument, it follows that
for each $\gamma \in \Gamma $, we have
$\gamma (\hat{K}(H_{1}))\subset \hat{K}(H_{2})\subset \hat{K}(H_{1})$, where the last inclusion holds since $H_1 \subset H_2$.  Thus since $\hat{K}(H_{1})$ is forward invariant under each $\gamma \in \Gamma$ and under each map in $H_1$, it is also forward invariant under each map in $H_2$.  We then conclude that
$\hat{K}(H_{1})\subset \hat{K}(H_{2})$, which together with the reverse inclusion already noted gives
$\hat{K}(H_{1})=\hat{K}(H_{2})$.
\end{proof}

\begin{definition}
Let $G$ be a rational semigroup and let $N$ be a positive integer.
We define $SH_N(G)$ to be the set of all $z\in \CC$ such that there
exists a neighborhood $U$ of $z$ such that for all $g \in G$ we have
$\deg(g:V \to U) \leq N$ for each connected component $V$ of
$g^{-1}(U)$.
\end{definition}

\begin{definition}
Let $G$ be a rational semigroup. We define $$UH(G)=\CC \setminus
\cup_{N=1}^\infty SH_N(G).$$
\end{definition}

\begin{remark}
For a rational semigroup $G$ we note that each $SH_N(G)$ is open and
thus $UH(G)$ is closed.
\end{remark}

\begin{remark}
For a rational semigroup $G$ we see that $UH(G) \subset P(G)$.  This
holds since for $z \notin P(G)$ and $U=B(z,\delta)$ such that $U
\cap P(G) =\emptyset$ it must be the case (by an application of the
Riemann-Hurwitz relation) that $\deg(g:V \to U)=1$ for each
connected component $V$ of $g^{-1}(U)$.
\end{remark}

\begin{remark}\label{cycles}
We note from Lemma 1.14 in \cite{Su2} that, the attracting cycles of
$g$, parabolic cycles of $g$, and the boundary of every Siegel disk
of $g$ are contained in $UH(\langle g\rangle )$, for any polynomial
$g$ with $\deg (g)\geq 2.$ Hence we may conclude that such points
are also in $UH(G)$ for any $G$ containing $g$.
\end{remark}

\begin{proof} [Proof of Theorem~\ref{semi-hyp}]
Assume the conditions stated in the hypotheses.  By the definition
of semi-hyperbolic, our goal is to show $J(G) \subset SH_K(G)$ for
some $K \in \N$.  We will show the equivalent statement that $J(G)
\cap UH(G) = \emptyset$. Since $UH(G) \subset P(G)$ and $P^*(G) \cap
J(G) \subset \Jmin(G)$, we have only to show $\Jmin(G) \subset \C
\setminus UH(G)$.

By Theorem~\ref{Jmin} or \cite{Su01,Su9} we know that $\sharp J_{\min }(G)\geq 3$. Thus
hypothesis (2) and Lemma~\ref{connect-to-Jf} imply
$\gamma^{-1}(\Jmin(G)) \cap \Jmin(G) = \emptyset$ for $\gamma \in
\Gamma$, which in turn implies (by Lemma~\ref{pullback-order})
$\gamma^{-1}(\Jmin(G))> \Jmin(G)$. Thus, for all $\gamma \in
\Gamma$,
\renewcommand{\theequation}{I}
\begin{equation}\label{star1}
\gamma^{-1}(J(G)) \cap A = \emptyset
\end{equation}
where $A=PH(\Jmin)$.

Since $\Gamma$ is compact in Poly, $d=\min_{\gamma \in \Gamma} dist
(\gamma^{-1}(J(G)), A)
>0$. By \eqref{star1} there exists $d_1>0$ such that for all $\gamma
\in \Gamma$, for all $z \in J(G)$, and all components $U$ of
$\gamma^{-1}(B(z,d_1))$ we have
\renewcommand{\theequation}{II}
\begin{equation}\label{star2}
U \cap B(A, d/2) = \emptyset.
\end{equation}
Now by Lemma~\ref{int=int} and by hypothesis (3) we have $UH(H) \cap
\C \subset P^*(H) \cap F(H) \subset \K(H) \cap F(H) =int \K(H) = int
\K(G) \subset F(G)$ and so, taking complements, $\Jmin(G) \subset \C
\setminus UH(H)$.

Claim:  There exists $b \in UH(H) \cap int\K(H)$.\\
Proof of claim:  Lemma~\ref{int=int} and Remark~\ref{facts} show
that $int \K(H) = int \K(G) \neq \emptyset$. Let $g_0 \in H$ and
consider the iterates $\{g_0^n\}$ at any $w \in int \K(H) \subset
F(H)$. Hypothesis (3) implies $UH(H) \cap \C \subset F(H)$ which
implies that $g_0$ cannot have a cycle of Siegel disks nor a
parabolic cycle (see Remark~\ref{cycles}). Thus by Sullivan's No
Wandering Domains Theorem the orbit $\{g_0^n(w)\}$ must be drawn
toward an attracting cycle in $\C$.  By replacing, if necessary,
$g_0$ by an iterate we may assume that $g_0^n(w)$ approaches a
finite fixed point $b$ of $g_0$.  Thus $b \in UH(H) \cap \C \subset
P^*(H) \cap F(H) \subset \K(H) \cap F(H)=int \K(H)$ which completes
the proof of the claim.

Now let $z \in \Jmin(G) \subset \C \setminus UH(H)$.  Then there
exists $\delta>0$ such that $B(z, 2 \delta) \subset \C \setminus
UH(H)$.  Since $g(UH(H)) \subset UH(H)$ for each $g \in H$, we must
have $g(b) \notin B(z,2\delta)$.  Since $H$ is normal at $b$, there
exists $\epsilon_1 >0$ such that $g \in H$ gives $g(B(b,
\epsilon_1)) \cap B(z, \delta) = \emptyset$, which implies
$g^{-1}(B(z, \delta)) \cap B(b, \epsilon_1) = \emptyset$.  Since $z
\in \C \setminus UH(H)$ there exists $\delta_1 < \delta$ and $N \in
\N$ such that for all $h \in H$ and for all components $V$ of
$h^{-1}(B(z, \delta_1))$ we have $\deg(h:V \to B(z, \delta_1)) \leq
N$.

Fix $h \in H$ and consider a component $V$ of $h^{-1}(B(z,
\delta_1))$ and note that the maximum principle implies that $V$ is
simply connected.  Let $\phi_{V,h}:B(0,1) \to V$ be the Riemann map
chosen such that $h \circ \phi_{V,h}(0)=z$.  By applying the
distortion Lemma~1.10 in~\cite{Su2}, there exists
$0<\delta_2<\delta_1$ such that the component $W$ of $(h\circ \phi
_{V,h})^{-1}(B(z,\delta _{2}))$ containing $0$ is such that $diam W
\leq c$ where $c>0$ is a small number independent of $h$, to be
specified later.

Note that, in the above, the set $V$ does depend on $h \in H$. Yet
for each $h \in H$, the set $\phi_{V,h}(B(0,1))=V$ does not meet
$B(b, \epsilon_1)$ and so the family $\{\phi_{V,h}\}_{h \in H}$ is
normal on $B(0,1)$. Thus
\renewcommand{\theequation}{III}
\begin{equation}\label{doubstar}
diam \phi_{V,h}(W) <d_1/10
\end{equation}
when $c$ is sufficiently small.

Let $g \in G$.  If $g \in H$, then (since $\delta_2< \delta_1$) we
have $\deg(g:V \to B(z, \delta_2)) \leq N$ where $V$ is any
component of $g^{-1}(B(z, \delta_2))$. If $g \notin H$, then we
write $g=h\gamma g_1$ where $g_1 \in G \cup \{id\}, h \in H \cup
\{id\}$ and $\gamma \in \Gamma$. Let $V_0$ be a component of
$\gamma^{-1}h^{-1}(B(z, \delta_2))$. Thus we have $\deg(h\gamma:V_0
\to B(z, \delta_2)) \leq NM$ where $M=\max_{\gamma \in \Gamma}\{\deg
\gamma \}$. By \eqref{doubstar} we have $diam \gamma (V_0) <
d_1/10$. By the definition of $d_1$ we have $V_0 \cap B(A,d/2) =
\emptyset$ and thus $V_0 \cap P(G) = \emptyset$. Using the maximum
principle applied to the polynomial $h\gamma $ implies $V_0$ is
simply connected and hence each branch of $g_1^{-1}$ is well defined
on $V_0$.  So for all components $V_1$ of $g^{-1}(B(z, \delta_2))$
we have $\deg(g:V_1 \to B(z, \delta_2)) \leq N M$.

In the above, $N$ depends on $z$, but what we have shown is that $z
\in \Jmin(G)$ implies $z \in \Jmin(G) \cap SH_N(H)$ for some $N$,
which in turn implies $z \in \Jmin(G) \cap SH_{NM}(G)$, thus giving
$z \notin UH(G)$. \end{proof}

\begin{proof}[Proof of Theorem~\ref{hyp}] The proof follows the same line
as the proof of Theorem~\ref{semi-hyp}.  We note that the usual
Koebe Distortion Theorem applies (without needing to invoke the
distortion Lemma 1.10 in~\cite{Su2}), and on the domains of interest
in the proof each $\gamma$ is one-to-one by hypothesis (4) and each
$h \in H$ is one-to-one by hypothesis (3).  We omit the details.
\end{proof}

\bibliographystyle{plain}
%

\begin{thebibliography}{10}

\bibitem{A}
Lars~V. Ahlfors.
\newblock {\em Conformal Invariants: Topics in Geometric Function Theory}.
\newblock McGraw-Hill, New York, 1973.

\bibitem{Be}
Alan~F. Beardon.
\newblock {\em Iterations of Rational Functions}.
\newblock Springer-Verlag, New York, 1991.

\bibitem{Bo-basin}
David~A. Boyd.
\newblock The immediate basin of attraction of infinity for polynomial
  semigroups of finite type.
\newblock {\em J. London Math. Soc. (2)}, 69(1):201--213, 2004.

\bibitem{Bruck}
Rainer Br{\"u}ck.
\newblock Geometric properties of {J}ulia sets of the composition of
  polynomials of the form {$z\sp 2+c\sb n$}.
\newblock {\em Pacific J. Math.}, 198(2):347--372, 2001.

\bibitem{BBR}
Rainer Br{\"u}ck, Matthias B{\"u}ger, and Stefan Reitz.
\newblock Random iterations of polynomials of the form {$z\sp 2+c\sb n$}:
  connectedness of {J}ulia sets.
\newblock {\em Ergodic Theory Dynam. Systems}, 19(5):1221--1231, 1999.

\bibitem{Bu1}
Matthias B{\"u}ger.
\newblock Self-similarity of {J}ulia sets of the composition of polynomials.
\newblock {\em Ergodic Theory Dynam. Systems}, 17(6):1289--1297, 1997.

\bibitem{Bu2}
Matthias B{\"u}ger.
\newblock On the composition of polynomials of the form {$z\sp 2+c\sb n$}.
\newblock {\em Math. Ann.}, 310(4):661--683, 1998.

\bibitem{CG}
Lennart Carleson and Theodore~W. Gamelin.
\newblock {\em Complex Dynamics}.
\newblock Springer-Verlag, New York, 1993.

\bibitem{FS}
John~Erik Forn{\ae}ss and Nessim Sibony.
\newblock Random iterations of rational functions.
\newblock {\em Ergodic Theory Dynam. Systems}, 11(4):687--708, 1991.

\bibitem{GR}
Z.~Gong and F.~Ren.
\newblock A random dynamical system formed by infinitely many functions.
\newblock {\em Journal of Fudan University}, 35:387--392, 1996.

\bibitem{HM1}
A.~Hinkkanen and G.J. Martin.
\newblock The dynamics of semigroups of rational functions {I}.
\newblock {\em Proc. London Math. Soc.}, 3:358--384, 1996.

\bibitem{N}
S.B.Nadler.
\newblock {\em Continuum Theory: An introduction}.
\newblock Marcel Dekker, 1992.

\bibitem{RS}
Rich Stankewitz.
\newblock {\em Completely invariant {J}ulia sets of rational semigroups}.
\newblock PhD thesis, University of Illinois, 1998.

\bibitem{RS1}
Rich Stankewitz.
\newblock Completely invariant {J}ulia sets of polynomial semigroups.
\newblock {\em Proc. Amer. Math. Soc.}, 127(10):2889--2898, 1999.

\bibitem{RS2}
Rich Stankewitz.
\newblock Completely invariant sets of normality for rational semigroups.
\newblock {\em Complex Variables Theory Appl.}, 40(3):199--210, 2000.

\bibitem{RS-Bobfest}
Rich Stankewitz.
\newblock Density of repelling fixed points in the Julia set of a Rational or Entire Semigroup.
\newblock {\em J. Difference Equ. Appl.}, (8 pages), to appear.

\bibitem{RS-JuliaSoftware}
Rich Stankewitz, W. Conatser, T. Butz, B. Dean, Y. Li, and K. Hart.
\newblock JULIA 2.0 Fractal Drawing Program.
\newblock http://rstankewitz.iweb.bsu.edu/JuliaHelp2.0/Julia.html
\bibitem{SSS}
Rich Stankewitz, Toshiyuki Sugawa, and Hiroki Sumi.
\newblock Some counterexamples in dynamics of rational semigroups.
\newblock {\em Ann. Acad. Sci. Fenn. Math.}, 29(2):357--366, 2004.

\bibitem{SS-AMC}
Rich Stankewitz and Hiroki Sumi.
\newblock Structure of {J}ulia sets of polynomial semigroups with
              bounded finite postcritical set.
\newblock {\em Appl. Math. Comput.}, 187(1):479--488, 2007.

\bibitem{Su1}
Hiroki Sumi.
\newblock On Hausdorff dimension of {J}ulia sets of hyperbolic rational
  semigroups.
\newblock {\em Kodai. Math.J.}, 21(1):10--28, 1998.

\bibitem{Su3}
Hiroki Sumi.
\newblock Skew product maps related to finitely generated rational semigroups.
\newblock {\em Nonlinearity}, 13:995--1019, 2000.

\bibitem{Su2}
Hiroki Sumi.
\newblock Dynamics of sub-hyperbolic and semi-hyperbolic rational semigroups
  and skew products.
\newblock {\em Ergod.Th.\& Dynam. Sys.}, 21:563--603, 2001.

\bibitem{Su7}
Hiroki Sumi.
\newblock Dimensions of {J}ulia sets of expanding rational semigroups.
\newblock {\em Kodai Mathematical Journal}, 28(2):390--422, 2005.
\newblock (See also http://arxiv.org/abs/math.DS/0405522.).

\bibitem{Su4}
Hiroki Sumi.
\newblock Semi-hyperbolic fibered rational maps and rational semigroups.
\newblock {\em Ergod.Th.\& Dynam. Sys.}, 26:893--922, 2006.

\bibitem{Su04} Hiroki Sumi.
\newblock Interaction cohomology of forward or backward self-similar systems.
\newblock {\em Adv. Math.}, 222 (2009) 729--781.

\bibitem{Su6}
Hiroki Sumi.
\newblock Dynamics of polynomial semigroups with bounded postcritical set in
  the plane.
\newblock {\em RIMS Kokyuroku}, 1447:198--215, 2005.
\newblock (Proceedings paper).

\bibitem{Su5}
Hiroki Sumi.
\newblock Dynamics of postcritically bounded polynomial semigroups and
  interaction cohomology.
\newblock {\em RIMS Kokyuroku}, 1447:227--238, 2005.
\newblock (Proceedings paper).


\bibitem{Su8}
Hiroki Sumi.
\newblock Random dynamics of polynomials and devil's staircase-like functions
  in the complex plane.
\newblock {\em Appl. Math. Comput.}, 187:489--500, 2007.
\newblock (Proceedings paper).

\bibitem{Su01}Hiroki Sumi.
\newblock Dynamics of postcritically bounded polynomial semigroups I: connected components
of the Julia sets.
\newblock  Discrete and Continuous Dynamical Systems Series A, 
Vol. 29, No. 3, 2011, 1205--1244.
\bibitem{Su02} Hiroki Sumi.
\newblock Dynamics of postcritically bounded polynomial semigroups II:
fiberwise dynamics and the Julia sets.
\newblock Preprint 2008, http://arxiv.org/abs/1007.0613.
\bibitem{Su03} Hiroki Sumi.
\newblock Dynamics of postcritically bounded polynomial semigroups III:
classification of semi-hyperbolic semigroups and random Julia sets
which are Jordan curves but not quasicircles.
\newblock Ergodic Theory Dynam. Systems (2010), 30, No. 6, 1869--1902.

\bibitem{Su05} Hiroki Sumi.
\newblock Random complex dynamics and semigroups of holomorphic maps.
\newblock Proc. London Math. Soc. (2011), 102 (1), 50--112. 
\bibitem{Su9}
Hiroki Sumi.
\newblock Dynamics of postcritically bounded polynomial semigroups.
\newblock Preprint 2006, http://arxiv.org/abs/math.DS/0703591.
\bibitem{SU2}
Hiroki Sumi and Mariusz Urba\'{n}ski.
\newblock The equilibrium states for semigroups of rational maps.
\newblock Monatsh. Math. 156 (2009), no. 4, 371--390.


\bibitem{SU1} Hiroki Sumi and Mariusz Urba\'{n}ski.
\newblock Real analyticity of Hausdorff dimension for
expanding rational semigroups.
\newblock  
Ergodic Theory Dynam. Systems (2010), Vol. 30, No. 2, 601-633. 
\bibitem{SU3} Hiroki Sumi and Mariusz Urba\'{n}ski.
\newblock Measures and dimensions of Julia sets of semi-hyperbolic
rational semigroups.
\newblock Discrete and Continuous Dynamical Systems Ser. A., 
Vol 30, No. 1, 2011, 313--363.

\bibitem{SY}
Y.Sun and C-C.Yang.
\newblock On the connectivity of the {J}ulia set of a finitely generated
  rational semigroup.
\newblock {\em Proc. Amer.Math.Soc.}, 130(1):49--52, 2001.

\bibitem{ZR}
W.~Zhou and F.~Ren.
\newblock The {J}ulia sets of the random iteration of rational functions.
\newblock {\em Chinese Sci. Bulletin}, 37(12):969--971, 1992.

\bibitem{Su11} 
Hiroki Sumi.  
\newblock Random complex dynamics and devil's coliseums. 
\newblock Preprint 2011, http://arxiv.org/abs/1104.3640
\end{thebibliography}

\end{document}